\documentclass{article}
\usepackage[utf8]{inputenc}

 \usepackage{tabu}
 
\usepackage{algorithm2e}
\usepackage{xcolor}
 
\SetAlFnt{\footnotesize}

\SetAlCapSty{xAlCapSty}

\SetCommentSty{xCommentSty}

\SetNlSty{mynlfont}{}{} 

\LinesNumbered

\SetSideCommentRight

\DontPrintSemicolon

\RestyleAlgo{algoruled}

\usepackage[driverfallback=dvipdfm,
bookmarks=true,    
bookmarksopen=true,
bookmarksopenlevel=\maxdimen,
bookmarksdepth=10
]{hyperref}
\hypersetup{
pdfcenterwindow=true, 
pdfpagelabels=true,
pagebackref=true,            
hypertexnames=true,
plainpages=false,
unicode=false,     
pdftoolbar=true,                
    pdfmenubar=true,          
    pdffitwindow=false,         
    pdfstartview={FitH},        
    pdftitle={Accuracy- Source Term Quadrature for 
 Third-Order Edge-Based Discretization},
    pdfauthor={Yi Liu and Hiroaki Nishikawa},    
    pdfsubject={Accuracy-Preserving Source Term Quadrature for 
 Third-Order Edge-Based Discretization},       
    pdfcreator={Hiroaki Nishikawa},   
    pdfproducer={Hiroaki Nishikawa}, 
    baseurl={http://www.hiroakinishikawa.com},
    pdfkeywords={}, 
    pdfnewwindow=true,     
    colorlinks=true,       
    linkcolor=black,       
    citecolor=black,       
    filecolor=black,       
    urlcolor=black,        
  breaklinks=true,
  hyperfigures=true,
  backref=true,
breaklinks=false, 
pdfpagelabels,      
pagebackref,        
hypertexnames=true, 
plainpages=false,   
naturalnames        
}
\usepackage[all]{hypcap}

\usepackage[thinlines]{easytable}

\usepackage{color}      
\usepackage{placeins}
\usepackage{indentfirst}
\usepackage{multirow}

\usepackage{amssymb}
\usepackage{amsmath,latexsym}
\usepackage{url}
\usepackage{enumerate}
\usepackage{graphicx}

\usepackage{booktabs}
 \newcommand{\ra}[1]{\renewcommand{\arraystretch}{#1}}
 
\usepackage{subcaption}
\usepackage{cleveref}
\usepackage{mathtools}

\usepackage{blindtext,graphicx}
\usepackage[absolute]{textpos}
\begin{document}


\title{\bf The QUICK Scheme is a Third-Order Finite-Volume Scheme with Point-Valued Numerical Solutions}

 \author{
{Hiroaki Nishikawa}\thanks{Associate Research Fellow ({hiro@nianet.org}),
 100 Exploration Way, Hampton, VA 23666 USA.}\\
  {\normalsize\itshape National Institute of Aerospace,
 Hampton, VA 23666, USA} 
}

\date{\today}
\maketitle

\begin{abstract} 
In this paper, we resolve the ever-present confusion over the QUICK scheme: it is a second-order scheme or a third-order scheme.  
The QUICK scheme, as proposed in the original reference [B. P. Leonard, Comput. Methods. Appl. Mech. Eng., 19, (1979), 59-98], is a third-order (not second-order) finite-volume scheme for the integral form of a general nonlinear conservation law with point-valued solutions stored at cell centers as numerical solutions. Third-order accuracy is proved by a careful and detailed truncation error analysis and demonstrated by a series of thorough numerical tests. 
The QUICK scheme requires a careful spatial discretization of a time derivative to preserve third-order accuracy for unsteady problems. Two 
techniques are discussed, including the QUICKEST scheme of Leonard. 
Discussions are given on how the QUICK scheme is mistakenly found to be second-order accurate. This paper is intended to serve as a reference to clarify any confusion about third-order accuracy of the QUICK scheme and also as the basis for clarifying economical high-order unstructured-grid schemes as we will discuss in a subsequent paper.
\end{abstract}

\section{Introduction}
\label{introduction}

This paper is a sequel to the previous paper \cite{Nishikawa_3rdMUSCL:2020}, where we discussed the third-order MUSCL scheme. 
The main motivation behind this work is to clarify some economical high-order unstructured-grid 
finite-volume schemes used in practical computational fluid dynamics solvers but largely confused in their mechanisms to achieve
third- and possibly higher-order accuracy (e.g., third-order U-MUSCL with $\kappa=1/2$ \cite{burg_umuscl:AIAA2005-4999}, $\kappa=1/3$ \cite{VLeer_Ultimate_III:JCP1977,VAN_LEER_MUSCL_AERODYNAMIC:J1985}, 
$\kappa=0$ \cite{katz_work:JCP2015,nishikawa_liu_source_quadrature:jcp2017}). 
In seeking the clarification, we have found that the confusion is rooted in the ever-present confusion over third-order convection schemes: the MUSCL scheme  \cite{VLeer_Ultimate_IV:JCP1977,VLeer_Ultimate_V:JCP1979} and the QUICK scheme \cite{Leonard_QUICK_CMAME1979}. To resolve the confusion, we started with the clarification of the third-order MUSCL 
scheme as given in Ref.\cite{Nishikawa_3rdMUSCL:2020}. In this paper, we will focus on the QUICK scheme.

The QUICK (Quadratic Upwind Interpolation for Convective Kinematics) scheme is a numerical scheme for convection equations originally developed by Leonard in 1979 \cite{Leonard_QUICK_CMAME1979}. 
Since then, it has become one of the most popular convection schemes widely used in practical computational fluid dynamics simulations and other related applications. The scheme is constructed based on a quadratic interpolation technique, the QUICK interpolation scheme, which is equivalent to the $\kappa$-reconstruction scheme \cite{VLeer_Ultimate_III:JCP1977,VAN_LEER_MUSCL_AERODYNAMIC:J1985} applied to point-valued solutions with $\kappa=1/2$. As such, there is no doubt that the QUICK interpolation scheme is third-order accurate since it is designed to be exact for quadratic functions. However, a controversy arose quickly and it still exists even today about the order of accuracy of the resulting convection scheme: third-order accurate
\cite{Yee_NASATM1985,Fletcher_v2_1991,HayaseHumphreyGreif:JCP1992,Darwish:NHT1993,Versteeg_CFD_Book1995,Chapter02:HandbookCFM1996,GjesdalTeigland:CNME1998,Roache_CFD_book,SongAmano:2001,Tkalich_QUCKEST_3rd:JH2006,Nacer_etal_HMT:2007,Zhang_etal:JCP2015,KumarLuoGasperOosterlee:JCP2018,DupuyToutantBataile:JCP2020} 
or second-order accurate 
\cite{Shyy_JCP1985,MackinonJohnson:IJNMF1991,Tsui:IJNMF1991,QUICK_is_2nd_ECCOMAS2006,Hirsch_VOL1_2nd_edition,Wesseling_CFD_Book,GhiasMittalDong:JCP2007,WatersonDeconinck_Convection:JCP2007,MerrickMalanRooyen:JCP2018}. Several authors including Leonard himself have attempted to resolve the confusion in the 1990's 
\cite{LeonardMokhtari:IJNMF1990,JohnsonMackinnon:CANM1992,ChenFalconer:AWR1994,Leonard_AMM1994,Leonard_AMM1995}, but the resolution does not seem to have
been achieved as we can find recent references stating that the QUICK scheme is second-order accurate 
\cite{QUICK_is_2nd_ECCOMAS2006,Hirsch_VOL1_2nd_edition,WatersonDeconinck_Convection:JCP2007}. 
This ever-present confusion seems to have been caused by the lack of rigorous accuracy studies: some references state that it is third-order but  
never verify the order of accuracy by numerical experiments \cite{Yee_NASATM1985,Fletcher_v2_1991,HayaseHumphreyGreif:JCP1992,Darwish:NHT1993,Versteeg_CFD_Book1995,Chapter02:HandbookCFM1996,GjesdalTeigland:CNME1998,Roache_CFD_book,SongAmano:2001,Tkalich_QUCKEST_3rd:JH2006,Nacer_etal_HMT:2007,Zhang_etal:JCP2015}; other references state that it is second-order by showing a truncation error analysis but often without taking into account the difference between a cell-averaged solution and a point-valued solution and again without numerical verifications \cite{MackinonJohnson:IJNMF1991,Tsui:IJNMF1991,QUICK_is_2nd_ECCOMAS2006,Hirsch_VOL1_2nd_edition}. 
 
The confusion would continue for a long time to come until a rigorous study is provided. 
In fact, the confusion has already been carried over to multi-dimensions and unstructured grids as mentioned at the beginning. In order to stop the spread of the confusion, we hereby provide a rigorous study consisting of a detailed truncation error analysis and a series of thorough numerical verification tests for steady and unsteady problems for a nonlinear conservation law with diffusion and forcing terms. It is the objective of this paper to clarify, once and for all, third-order accuracy of the QUICK scheme, and to establish the foundation for clarifying third-order unstructured-grid schemes as we will discuss in a subsequent paper. 

The truth is that the QUICK scheme is {\it a third-order finite-volume discretization of the integral form of a conservation law with point-valued solutions stored as numerical solutions at cell centers}, which is exactly how it was introduced in Ref.\cite{Leonard_QUICK_CMAME1979} and clearly stated again later in Ref.\cite{Leonard_AMM1994}. Therefore, it differs from the third-order MUSCL scheme only by the definition of the numerical solution: the QUICK scheme with point-valued solutions and the MUSCL scheme with cell-averaged solutions \cite{Nishikawa_3rdMUSCL:2020}. The lack of distinction between the point value and the cell average is 
one of the major sources of confusion. Another major source of confusion is the lack of distinction between a point-valued operator and a cell-averaged operator: 
 the QUICK and MUSCL schemes both approximate the latter. Yet another is the lack of understanding that the QUICK scheme needs a careful spatial discretization of the time derivative term to achieve third-order accuracy in space, which was already pointed out in the original 1979 paper \cite{Leonard_QUICK_CMAME1979} but has often been missed. 
 
Leonard has correctly recognized the distinction between a point-valued operator and a cell-averaged operator, and repeatedly pointed it out in subsequent papers \cite{LeonardMokhtari:IJNMF1990,Leonard_AMM1994,Leonard_AMM1995}. However, in his analyses, the numerical solution had always been considered as a point value. It seems to be the main reason that the $\kappa=1/3$ MUSCL scheme did not look to him as a third-order scheme until he published Ref.\cite{Leonard_AMM1994}, where he considered cell-averaged solutions as well and concluded that the QUICK scheme was third-order with point values and the MUSCL scheme is third-order with cell averages. He also showed that the $\kappa=1/3$ scheme (with point-valued solutions) was a third-order finite-difference scheme in Ref.\cite{Leonard_AMM1994}, but it is actually true only for linear equations on one-dimensional uniform grids \cite{Shu_WENO_SiamReview:2009}. This $\kappa=1/3$ scheme appears also in Ref.\cite{JohnsonMackinnon:CANM1992}, where it is called QUICK-FD. 
Ref.\cite{Wesseling_CFD_Book} analyzes the MUSCL scheme with Van Leer's $\kappa$-reconstruction and concludes that the QUICK scheme corresponds to $\kappa=1/2$ and is only second-order accurate. It also discusses the QUICKEST scheme and associates it with the $\kappa=1/3$ scheme, which are indeed related \cite{LeonardMokhtari:IJNMF1990}, but it never discusses the finite-volume method with point-valued solutions, for which the QUICK scheme is third-order. 
It is noted that all these analyses have been performed for a linear equation; no analyses are found in the literature, to the best of the author's knowledge, of the QUICK scheme for a nonlinear equation, which can generate an additional confusion as we will discuss later. In this paper, we will focus on a general nonlinear conservation law.

The need for a spatial discretization of the time-derivative discretization is evident from the fact that 
the QUICK scheme is based on the integral form having the time derivative of the cell-averaged solution. Hence, it does not
immediately provide an evolution equation for a point-valued numerical solution, and the time derivative of the cell-averaged solution must be expressed in terms of point-valued solutions even before it is discretized in time. Leonard already recognized the need for such a consistent treatment in 1979 \cite{Leonard_QUICK_CMAME1979}, and developed a third-order unsteady scheme, the QUICKEST (QUICK with Estimated Streaming Terms) scheme, which incorporates a consistent time-derivative treatment in such a way to create an explicit convection scheme. This consistent treatment, however, has not been well understood among those who claim that the QUICK scheme is second-order accurate, as mentioned in Ref.\cite{LeonardMokhtari:IJNMF1990}. In this paper, we generalize the QUICKEST scheme as a semi-discrete scheme and show that it can be integrated in time by any time-stepping scheme. However, we also show that it requires a flux reconstruction to achieve third-order accuracy for nonlinear equations and furthermore that third-order can be achieved only for uniform grids. As a more general technique to preserve third-order accuracy for unsteady problems, we will consider a coupled mass-matrix formulation. It can be thought of as a method derived from a general deconvolution approach \cite{Denaro:IJNMF1996,Denaro_CMAME:2015}, which extends systematically to arbitrarily high-order and has been demonstrated for unstructured triangular grids also \cite{Denaro:IJNMF2002}. 

Yet another confusion exists about the accuracy of the QUICK scheme at a face. 
Ref.\cite{Hirsch_VOL1_2nd_edition} states that the QUICK scheme will be third-order accurate if the solution values are stored at faces rather than at cells (or nodes), which is consistent with the claim that the QUICK scheme has a third-order truncation error at a 
face \cite{LeonardMokhtari:IJNMF1990,ChenFalconer:AWR1994}. This claim, to the best of the author's knowledge, has never been verified numerically nor correctly proved analytically. Note also that the scheme is incomplete since it is not clear how to advance the solutions at cell centers (or at nodes) in time, from 
which the QUICK interpolation is performed at a face. Ref.\cite{LeonardMokhtari:IJNMF1990} presents an analysis and concludes that the solutions interpolated at the left and right faces are both third-order accurate and so the convection scheme based on their difference is third-order accurate. However, it does not correctly predict the order of accuracy because the solutions at the left and right faces are expanded independently around the left and right face points, respectively. The Taylor expansion must be performed at a single point, where the scheme is defined, to correctly derive its truncation error, as pointed out later in Ref.\cite{JohnsonMackinnon:CANM1992}. In this case, the point should be the right face, at which third-order accuracy is claimed. However, it is easy to verify by Taylor expansions around the right face that the QUICK interpolated solutions at the left and right faces are first- and third-order accurate, respectively. Therefore, a convection scheme defined at the right face based on their difference is only first-order accurate. In this paper, we do not discuss such a face-centered scheme because it is not complete. Interested readers are referred to 
Refs.\cite{VLeer_Ultimate_IV:JCP1977,eymann_roe:AIAA2011,nishikawa_roe_active_flux_advdiff:CF2016} for a high-order scheme incorporating an upwind evolution of a flux stored at a face, which could be 
considered as an example of a complete face-based scheme.

 The paper is organized as follows.
In Section 2, we will describe a target conservation law and its integral form, and clarify the fundamental difference between the MUSCL and QUICK schemes.
In Section 3, we provide a detailed description of the QUICK scheme designed for steady problems and two approaches to preserving third-order accuracy for unsteady problems, including the QUICKEST scheme.  
In Section 4, we provide a detailed truncation error analysis and prove that the QUICK and its unsteady versions are all third-order accurate.
In Section 5, we discuss the truncation error of the QUICK scheme as a finite-difference scheme approximating the differential form of a conservation law instead of the integral form, and discuss how confusing it is.
In Section 6, we provide a summary of various QUICK schemes and their formal orders of accuracy.
In Section 7, numerical results are presented.
In Section 8, concluding remarks are given.

\section{Target Equation and Exact Integral Form}
\label{QUICK_I_integral_form}
 
 Consider a one-dimensional conservation law, including a diffusion term and a forcing term:
\begin{eqnarray}
u_t  + f_x = \nu u_{xx} + {s(x)},
\label{diff_form}
\end{eqnarray}
where ${u}$ is a solution variable, $f=f(u)$ is a nonlinear flux, $\nu$ is a positive constant diffusive coefficient, ${s}(x)$ is a forcing term, and 
the subscripts $t$ and $x$ denote the partial derivatives with respect to time and space, respectively. 
The diffusion and forcing terms are included because the treatment of these terms is very important to the QUICK scheme: third-order 
accuracy can be easily lost with a mistreatment of these terms. 

  \begin{figure}[t]
    \centering
          \begin{subfigure}[t]{0.485\textwidth}
        \includegraphics[width=0.99\textwidth]{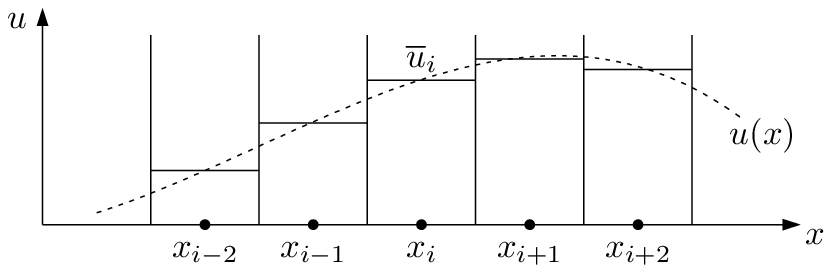} 
          \caption{Cell-averaged solutions: MUSCL.} 
          \label{fig:oned_fv_data_p0}%
      \end{subfigure}
          \begin{subfigure}[t]{0.485\textwidth}
        \includegraphics[width=0.99\textwidth]{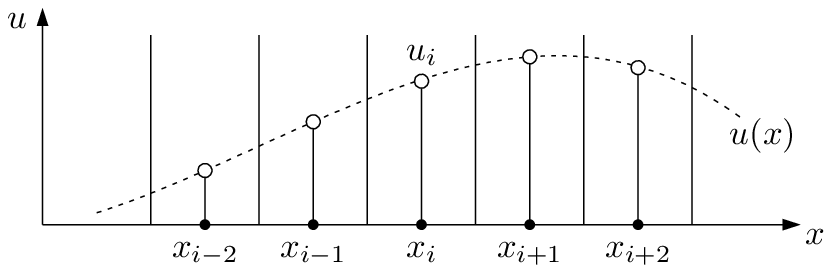} 
          \caption{ Point-value solutions: QUICK.}
        \label{fig:oned_fv_data_point}%
      \end{subfigure}
            \caption{
\label{fig:solution_typ_on_oned_uniform}%
Cell-averaged solutions at cells and point-valued solutions at cell centers for a smooth function $u(x)$ on a uniform grid.
} 
\end{figure}

As described in Ref.\cite{Leonard_QUICK_CMAME1979}, the QUICK scheme is constructed based on the integral form of a conservation law.
Consider a one-dimensional grid with a uniform spacing $h$: $x_{i+1} -x_i = h$, $i=1,2,3, \cdots, n$, where $n$ is an integer, as shown in 
Figure \ref{fig:oned_fv_data_p0}. We integrate the differential form (\ref{diff_form}) over a control volume around $x=x_i$,
$x \in [x_{i-1/2}, x_{i+1/2} ] =  [x_{i} - h/2, x_{i}+h/2 ] $, and obtain
\begin{eqnarray}
 \frac{d  \overline{u}_i }{dt}  +    \frac{1}{h} \int_{x_i - h/2}^{x_i+h/2} f_x \, dx    =  
  \frac{1}{h} \int_{x_i - h/2}^{x_i+h/2} \nu u_{xx}  \, dx  +  \overline{s}_i,
\end{eqnarray}
or
\begin{eqnarray}
 \frac{d  \overline{u}_i }{dt}     +  \frac{1}{h}  \left[ \left\{ f(u_{i+1/2}) -\nu  (u_x)_{i+1/2} \right\} - 
 \left\{  f(u_{i-1/2})  -\nu  (u_x)_{i-1/2}  \right\}  \right]   =  \overline{s}_i,  
 \label{integral_form_02}
\end{eqnarray}
where $\overline{u}$ and $\overline{s}$ are cell-averaged solution and forcing term:
\begin{eqnarray}
  \overline{u}_i   =   \frac{1}{h}  \int_{x_i-h/2}^{x_i+h/2} {u} \, dx ,\quad
  \overline{s}_i   =   \frac{1}{h} \int_{x_i-h/2}^{x_i+h/2} {s}(x) \, dx , 
\end{eqnarray}
and $u_{i-1/2}$ and $u_{i+1/2}$ are the point-valued solutions at the left and right faces, respectively.
Note that the integral form is exact and no approximation has yet been made.


An natural choice for the numerical solution is the cell-average $\overline{u}_i$, which  
leads to the MUSCL scheme as discussed in Ref.\cite{Nishikawa_3rdMUSCL:2020}. But the numerical solution
does not have to be the cell-average. In principle, we can choose anything we want and can still discretize the integral form.
In Ref.\cite{Leonard_QUICK_CMAME1979}, the point-valued solution was chosen as the numerical solution and stored at cell centers, and then 
the QUICK scheme was proposed as a finite-volume discretization of the integral form. See Figure \ref{fig:oned_fv_data_point}.  
As one can easily see, the point-valued solution needed at a face for the flux evaluation is now interpolated from the point-valued solutions stored at cells, 
instead of reconstructed from cell averages as in the MUSCL scheme. This is exactly what the QUICK interpolation scheme does 
for a quadratic solution, thus leading to a third-order convection scheme. Therefore, the only difference between the QUICK ands MUSCL schemes
is the type of the numerical solution. See Table \ref{Tab.MUSCL_QUICK}. 
The time derivative remains a cell-average, 
as clearly seen in Equation (\ref{integral_form_02}), but the numerical solution is a point value; this is the reason that a careful discretization is needed 
for unsteady problems, as we will discuss later.

To derive a truncation error, we must be clear about the fact that the target equation approximated by the QUICK scheme is the integral form: 
\begin{eqnarray}
  \frac{d  \overline{u}_i }{dt}  + (\overline{f_x})_i  -(\overline{\nu u_{xx}})_i  =   \overline{s}_i, 
  \label{integral_form_operator}
\end{eqnarray}
where
\begin{eqnarray}
(\overline{f_x})_i  =    \frac{1}{h} \int_{x_i - h/2}^{x_i+h/2} f_x \, dx,  \quad 
(\overline{\nu u_{xx}})_i =   \frac{1}{h} \int_{x_i - h/2}^{x_i+h/2} \nu u_{xx} \, dx,
\end{eqnarray}
and expand the QUICK scheme to find how well it approximates this equation, not the differential form (\ref{diff_form}). 
As we will show, the truncation error is third-order and thus the QUICK scheme is third-order accurate. There is no confusion. 
A confusion arises when the QUICK scheme is mistakenly considered as a finite-difference scheme as we will discuss later.

\begin{table}[t]
\ra{1.25}
\begin{center}
\begin{tabular}{rccccc}\hline\hline 
\multicolumn{1}{r}{ }                               &      
\multicolumn{1}{r}{\small Numerical solution}     &
\multicolumn{1}{r}{\small Target equation}    &
\multicolumn{1}{r}{\small Flux at a face}  &
\multicolumn{1}{r}{\small Time derivative}   &
\multicolumn{1}{r}{\small Forcing/source term} 
  \\ \hline  
\small MUSCL& \small    Cell average  &\small   Integral form  & \small  Point value&   \small  Cell average &   \small  Cell average   \\ 
\small QUICK  & \small  {\bf  \color{red}  Point value}  &  \small Integral form  &\small   Point value&  \small   Cell average & \small    Cell average 
  \\  \hline  \hline  
\end{tabular}
\caption{Comparison of MUSCL and QUICK schemes.
}
\label{Tab.MUSCL_QUICK}
\end{center}
\end{table}

\section{Third-Order QUICK Scheme}
\label{third_order_QUICK_I}
 
In this section, we follow the original paper \cite{Leonard_QUICK_CMAME1979} and describe the QUICK scheme
as a finite-volume discretization constructed from point-valued solutions stored at cell centers as numerical solutions.
 To preserve third-order accuracy for unsteady problems, we need to carefully
discretize the time derivative term. Two approaches are discussed: the QUICKEST scheme of Leonard \cite{Leonard_QUICK_CMAME1979} and
a coupled mass-matrix formulation.

\subsection{QUICK: Third-order spatial discretization}
\label{third_order_QUICK_space}

\begin{figure}[t]
    \centering 
        \includegraphics[width=0.58\textwidth]{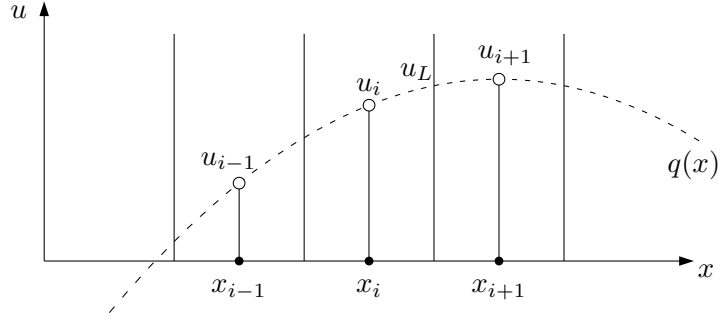} 
            \caption{
\label{fig:oned_fv_data_quad_inter}%
Quadratic interpolation over the three cell centers to obtain the solution on the left side of the right face $i+1/2$.
} 
\end{figure}

Let us denote a point-valued numerical solution at a cell center $x=x_i$ by $u_i$, and consider the following finite-volume discretization: 
\begin{eqnarray}
  \frac{d  \overline{u}_i }{dt}  +   \frac{1}{h}  \left[  F_{i+1/2} -F_{i-1/2}  \right]  =  \overline{s}_i,
\label{fv_exact_form_discretized}
\end{eqnarray}
where $F_{i\pm1/2} = F(u_{i \pm /2,L},(u_x)_{ \pm 1/2,L},u_{i \pm 1/2,R},(u_x)_{i \pm /2,R} )$, and $F$ denotes a numerical flux as a function of two point-valued solutions and derivatives interpolated at a face. In this paper, we consider the following numerical
flux, e.g., at the right face with the subscript $i+1/2$ dropped, 
\begin{eqnarray}
F(u_L,u_R)  =  F^{c}(u_L,u_R)   + F^d(u_L, (u_x)_L, u_R, (u_x)_R),
\label{upwnd_flux}
\end{eqnarray}
where $F^{c}$ is an upwind convective flux,
\begin{eqnarray}
F^{c}(u_L,u_R)  = \frac{1}{2} \left[  f(u_L)    + f(u_R)     \right]  - \frac{D}{2}  ( u_R - u_L), 
\label{upwnd_flux_convection}
\end{eqnarray}
with the dissipation coefficient $D= |\partial f/ \partial u|$, $F^{d}$ is the alpha-damping diffusive flux \cite{nishikawa:AIAA2010,nishikawa_general_principle:CF2011},
\begin{eqnarray}
F^d(u_L,u_R)  = -  \frac{1}{2} \left[   \nu (u_x)_L   +  \nu (u_x)_R  \right]  + \frac{\nu \alpha}{2 h}  ( u_R - u_L),
\label{upwnd_flux_diff}
\end{eqnarray}
with $\alpha$ as a constant damping coefficient to be determined later, and $u_L$, $ (u_x)_L$, $u_R$, and $(u_x)_R$ are interpolated point-valued solutions and consistently evaluated derivatives at a face from the left and right  cells, respectively. 
 For the solution interpolation, we consider Van Leer's $\kappa$-reconstruction scheme \cite{VAN_LEER_MUSCL_AERODYNAMIC:J1985,VLeer_Ultimate_III:JCP1977} applied to point-valued solutions:
\begin{eqnarray}
u_L &=& \frac{1}{2}  \left(  {u}_{i} +   {u}_{i+1}   \right)- \frac{1-\kappa}{4} \left(   {u}_{i+1}  -2 {u}_{i}+  {u}_{i-1} \right),  
\label{kappa_uL_simple} \\ [2ex]
u_R &=& \frac{1}{2}  \left(  {u}_{i+1} +  {u}_{i}   \right)- \frac{1-\kappa}{4} \left(  {u}_{i+2}  -2  {u}_{i+1}+ {u}_{i} \right).
\label{kappa_uR_simple}
\end{eqnarray}
Our focus is on the QUICK interpolation scheme of Leonard, corresponding to $\kappa=1/2$:
\begin{eqnarray}
u_L =
\frac{1}{2}  \left( {u}_{i} +   {u}_{i+1}   \right)- \frac{1}{8} \left(  {u}_{i+1}  -2 {u}_{i}+  {u}_{i-1} \right), \label{QUICK_uL} \\ [1.5ex]
u_R =   \frac{1}{2}  \left( {u}_{i+1} +   {u}_{i}   \right)- \frac{1}{8} \left(  {u}_{i+2}  -2 {u}_{i+1}+  {u}_{i} \right), 
\label{QUICK_uR} 
\end{eqnarray}
each of which quadratically interpolates the point-valued solution, e.g., see Figure \ref{fig:oned_fv_data_quad_inter} for $u_L$.
In what follows, we will keep the general form as in Equations (\ref{kappa_uL_simple}) and (\ref{kappa_uR_simple}) because the choice $\kappa=1/3$ is also relevant to the QUICK scheme as will be discussed later.
 \newline
 \newline
\noindent {\bf Remark}: A similar reconstruction scheme is used in Ref.\cite{yang_harris:AIAAJ2016}, which adds
extra terms to the $\kappa$-reconstruction scheme with an additional parameter. But when the additional parameter is zero, it does not reduce to the 
$\kappa$-reconstruction scheme: 
\begin{eqnarray}
u_L = \frac{1}{2}  \left(  {u}_{i} +   {u}_{i+1}   \right)- \frac{1-\kappa}{4} \left(   {u}_{i+1}  -2 {u}_{i}+  {u}_{i-1} \right)
+ \frac{1}{32}  \left(   {u}_{i+2}  -2 {u}_{i}+  {u}_{i-2} \right),
\end{eqnarray}
which is claimed to achieve third-order accuracy for a linear convection equation with $\kappa=-1/6$ in the point-valued solution.
Hence, this interpolation scheme is built upon a slightly different interpolation scheme; it involves additional solution values $ {u}_{i-2}$ and ${u}_{i+2}$, and therefore
will result in a larger stencil than that of the QUICK scheme considered here. This particular scheme is beyond the scope of this paper, and will not be 
discussed in the rest of the paper.
 \newline
 \newline
At each face, we also need to compute the derivatives for the diffusive flux. For third-order accuracy on uniform grids, we wish to derive these derivatives
exactly for a quadratic function. In consistent with the QUICK interpolation, the derivative at a face must be obtained from 
a quadratic interpolating polynomial. Consider the following polynomials,
\begin{eqnarray}
q_i(x) &=& {u}_{i} + (u_x)_{i}  (x - x_{i} ) + \frac{1}{2} (u_{xx})_{i}  (x - x_{i} )^2 , \\ [2ex]
q_{i+1}(x)  &=& {u}_{i+1} +  (u_x)_{i+1} (x-x_i) +  \frac{1}{2}(u_{xx})_{i+1} (x - x_{i+1} )^2,
\label{quick_01}
\end{eqnarray}
where 
\begin{eqnarray} 
 (u_x)_{i}  = \frac{  {u}_{i+1} -  {u}_{i-1}  }{2h}, \quad 
 (u_{xx})_{i}  =  \frac{  {u}_{i+1}  -2 {u}_{i}+  {u}_{i-1}  }{h^2}, \quad
 (u_x)_{i+1}  = \frac{  {u}_{i+2} -  {u}_{i}  }{2h}, \quad 
 (u_{xx})_{i+1}  =  \frac{  {u}_{i+2}  -2 {u}_{i+1}+  {u}_{i}  }{h^2},
 \label{fd_aprpox_at_ip1}
\end{eqnarray}
each of which quadratically interpolates three point-valued solutions: e.g., $q_i(x_{i-1}) = u_{i-1}$, $q_i(x_i) = u_i$, $q_i(x_{i+1}) = u_{i+1}$.
Then, we derive the derivatives $(u_x)_L$ and $ (u_x)_R$ by differentiating these polynomials and evaluating them at the face:
\begin{eqnarray}
(u_x)_L  = \left. \frac{d q_i(x)}{d x} \right|_{x=x_i+h/2} 
                 = \frac{  {u}_{i+1} -  {u}_{i}  }{h} , 
                 \quad
(u_x)_R  = \left. \frac{d q_{i+1}(x)}{d x} \right|_{x=x_{i+1}-h/2}  
                 =  \frac{  {u}_{i+1} -  {u}_{i}  }{h},
\end{eqnarray}
indicating that the derivative of the quadratic interpolation polynomial is continuous across the face.
Using these derivatives in the flux (\ref{upwnd_flux_diff}), we obtain the following diffusion scheme:
\begin{eqnarray}
    \frac{1}{h}  \left[  F^d_{i+1/2} -F^d_{i-1/2}  \right]  &=&  -  \frac{ \nu }{8 h^2} 
 \left[ 
    \alpha(\kappa-1) u_{i-2}  - 4 \left\{     \alpha (\kappa-1) - 2 \right\}  u_{i-1} 
  +  2 \left\{   3 \alpha (\kappa-1) - 8   \right\} u_i  \right. \nonumber  \\ [1ex]
 &-&
 \left.
4 \left\{    \alpha (\kappa-1)  - 2  \right\} u_{i+1}   +   \alpha(\kappa-1) u_{i+2}
    \right],
\label{fv_exact_form_discretized_diff_kappa_alpha}
\end{eqnarray}
which is a two-parameter-family scheme with $\alpha$ and $\kappa$ as the parameters. Note that we have assumed that the same $u_L$ and $u_R$ are used for
the dissipation term in the convective flux and the damping term in the diffusive flux; they can be computed, if one wishes, with a different value of $\kappa$ for the diffusive flux, to make the diffusion scheme independent of the choice of the convective scheme. 
To achieve third-order accuracy (fourth-order accuracy as a diffusion scheme) in uniform grids, however, there is 
a relationship that must be satisfied by $\alpha$ and $\kappa$ as we will show later. 
 
\subsection{Consistent spatial discretization of time derivatives}
\label{third_order_QUICKIEST}

\subsubsection{The need for spatial discretization of the time derivative}

At this point, the scheme (\ref{fv_exact_form_discretized}) is still in the form
of the cell-average evolution equation. Therefore, it cannot be used directly to advance the point-valued solution $u_i$ in time.
To derive a point-valued evolution equation, we consider the relationship between the the cell average and the point value:
\begin{eqnarray}
 \overline{u}_{i }  
 =  u_i + \frac{h^2}{24} (u_{xx})   + O(h^4),
 \label{ca_exact_expansion}
\end{eqnarray}
where $(u_{xx})$ denotes the second-derivative at $x=x_i$, 
which can be easily obtained by cell-averaging the Taylor expansion of a smooth point-valued solution \cite{Nishikawa_3rdMUSCL:2020}. 
Then, we write the time derivative of the cell-averaged solution as 
\begin{eqnarray}
 \frac{d   \overline{u}_i }{dt} 
 =    \frac{d   {u}_i }{dt} + \frac{h^2}{24} \left(   \frac{d  (u_{xx})  }{d t}  \right)    + O(h^4),
 \label{ca_exact_expansion_dudt}
\end{eqnarray}
and write the scheme (\ref{fv_exact_form_discretized}) as
\begin{eqnarray}
  \frac{d   {u}_i }{dt}  + \frac{h^2}{24}  \left(   \frac{d  (u_{xx})  }{d t}  \right)     + \frac{1}{h}  [ F_{i+1/2}  - F_{i-1/2}  ]   = \overline{s}_i,
      \label{pointwise_scheme0000}
\end{eqnarray}
or
\begin{eqnarray}
  \frac{d   {u}_i }{dt}  + \frac{h^2}{24} \left(   \frac{d   {u} }{dt}  \right)_{\!\! xx}    + \frac{1}{h}  [ F_{i+1/2}  - F_{i-1/2}  ]   = \overline{s}_i,
      \label{pointwise_scheme0001999}
\end{eqnarray}
which is accurate enough to update the point-valued solution ${u}_i$ with third-order accuracy.  
The conversion of the cell average to the point value is often called the deconvolution \cite{Denaro:IJNMF1996,Harten:jcp1987}. The above deconvolution 
is a simple example of a general deconvolution approach described in Refs.\cite{Denaro:IJNMF1996,Denaro_CMAME:2015,Denaro:IJNMF2002}: the operator 
$1 + h^2 \partial_{xx}/24$ considered as applied to $ d u_i /dt $ in the above equation is a truncated deconvolution operator of fourth-order 
accuracy \cite{StefanoDenaroRiccardi:IJNMF2001}. 
Here, to achieve third-order accuracy, we need to accurately discretize the second-derivative of the time derivative $ \left(   \frac{d   {u} }{dt}  \right)_{\!\! xx}  $  in space 
or discretize the second-derivative $ (u_{xx})$ with sufficient accuracy and then take its time derivative. 
In this paper, we consider two approaches: the QUICKEST scheme of Leonard and a coupled mass-matrix scheme.

\subsubsection{QUICKEST}
\label{quickest_explained}
 
Leonard already recognized the need for a consistent discretization of the time-derivative term and 
developed a fully discrete explicit time-accurate QUICK scheme  
in his 1979 paper \cite{Leonard_QUICK_CMAME1979}, which is called the QUICKEST scheme. 
Here, we follow his derivation for a linear convective equation $u_t+f_x =0$, $f=au$ with a global constant $a$, but keep the semi-discrete form in order to be
able to independently apply a high-order time-integration scheme. The basic idea for deriving an explicit scheme is to convert
the curvature term in Equation (\ref{pointwise_scheme0001999}) to a purely spatial derivative term by using $u_t+f_x=0$,
\begin{eqnarray} 
 \frac{h^2}{24} \left(   \frac{d   {u} }{dt}  \right)_{\!\! xx} 
=
-   \frac{h^2}{24} \left(   f_x  \right)_{xx}
= -  \frac{h^2}{24}  a  ( u_{xx}  )_x,
\label{QUICKEST_curv_conversion}
 \end{eqnarray}
 and approximate it as 
 \begin{eqnarray}
 \frac{h^2}{24} \left(   \frac{d   {u} }{dt}  \right)_{\!\! xx}   \approx
-   \frac{a}{24 h}  \left[   (  u_{i+1} -2 u_i + u_{i-1} ) - ( u_{i} -2 u_{i-1} + u_{i-2} )      \right].
\label{low_order_curvature_diff_QUICKEST}
 \end{eqnarray}
Note that this is a first-order accurate approximation at $x=x_i$ (which can be easily shown by a Taylor expansion around $x=x_i$), but
the resulting scheme can be third-order. Substituting it into Equation (\ref{pointwise_scheme0001999}), we arrive at the following semi-discrete equation,
\begin{eqnarray}
  \frac{d   {u}_i }{dt}  +  \frac{1}{h}  [ \tilde{F}^c_{i+1/2}  - \tilde{F}^c_{i-1/2}  ]   =  0,
      \label{pointwise_scheme1a}
\end{eqnarray}
where the modified flux $\tilde{F}$ is given, for a pure upwind flux ${F}^c = a u_L$ with $\kappa=1/2$, by
\begin{eqnarray}
\tilde{F}^{c}(u_L,u_R)  
& =&  a u_L -  \frac{a}{24} (  u_{i+1} -2 u_i + u_{i-1} )  \\ [1ex]
&=& a \left[  \frac{1}{2}  \left( {u}_{i} +   {u}_{i+1}   \right) - \frac{1}{8} \left(  {u}_{i+1}  -2 {u}_{i}+  {u}_{i-1} \right)   
-   \frac{1}{24} (  u_{i+1} -2 u_i + u_{i-1} ) \right]  \\ [1ex]
&=& 
 a \left[  \frac{1}{2}  \left( {u}_{i} +   {u}_{i+1}   \right)  -   \frac{1}{6} (  u_{i+1} -2 u_i + u_{i-1} ) \right].
\end{eqnarray}
This is the QUICKEST scheme of Leonard in the semi-discrete form. 
Notice that this corresponds to the original flux ${F}^{c}$ with $\kappa=1/3$, and thus the QUICKEST scheme can be written as
\begin{eqnarray}
  \frac{d   {u}_i }{dt}  +  \frac{1}{h}  [  {F}^c_{i+1/2}  - {F}^c_{i-1/2}  ]   =  0, \quad \kappa = \frac{1}{3},
      \label{pointwise_scheme1aaaa}
\end{eqnarray}
which is a convenient form directly applicable to nonlinear equations. This third-order scheme is equivalent to 
a third-order finite-difference scheme ($\kappa=1/3$) derived by Van Leer in 
1977 \cite{VLeer_Ultimate_III:JCP1977} (the $\kappa$-reconstruction scheme was
first developed for a finite-difference scheme). For this reason, the QUICKEST scheme is sometimes said to
be equivalent to or is rediscovered from the third-order MUSCL 
scheme ($\kappa=1/3$) \cite{Wesseling_CFD_Book,LeonardMokhtari:IJNMF1990,Leonard_AMM1994}; but such is not true because the MUSCL scheme is based on
 cell-averaged solutions while the QUICKEST scheme is based on point-valued solutions. 
 A fully discrete scheme constructed by Leonard in Ref.\cite{Leonard_QUICK_CMAME1979} can be obtained
by integrating the time-derivative term in time and 
quadratically extrapolating the fluxes to the next time level at both faces (see Ref.\cite{Leonard_QUICK_CMAME1979}).
A similar formulation can be found in Ref.\cite{DavisMoore:JFM1982}, which extends the QUICKEST scheme to two dimensions. 
In this paper, we keep the semi-discrete form and show that the core idea of the QUICKEST scheme is more general and allows the use of any time-integration scheme.  

The QUICKEST scheme is a finite-difference scheme 
as implied by the point-valued time-derivative in Equation (\ref{pointwise_scheme1aaaa}), thus the flux balance approximating the flux derivative at the cell center. 
As such, the forcing term (if present) must be evaluated at the cell center, not cell averaged, and any other term will have to be discretized as a point-valued 
approximation. This sudden switch from finite-volume to finite-difference has been caused by the
low-order approximation (\ref{low_order_curvature_diff_QUICKEST}). If the conversion (\ref{QUICKEST_curv_conversion}) is approximated more accurately, 
then we would keep third-order accuracy as a finite-volume scheme. For example, the general deconvolution approach of Denaro \cite{Denaro:IJNMF1996} can be used to directly derive high-order explicit schemes for a general nonlinear conservation law. Or the second-derivative of the point time-derivative in Equation (\ref{pointwise_scheme0001999}) can be discretized by a central difference formula applied to the cell-average time-derivatives replaced by the spatial residuals. This explicit scheme is due to Van Leer and can be used to generate high-order explicit schemes. It can be easily implemented for a general nonlinear conservation law; we will provide a brief description and present numerical results later in Section \ref{subsec:unsteady}. 
Typically, such explicit high-order schemes will enlarge the residual stencil, beyond the typical five-point stencil, because it involves high-order derivatives of the spatial operator. In this paper, we will focus on the schemes defined within the five-point stencil.

It should be noted that the QUICKEST scheme (\ref{pointwise_scheme1aaaa}) is indeed third-order accurate for linear equations \cite{Wesseling_CFD_Book,Leonard_AMM1994} (on uniform grids), but only second-order accurate for nonlinear equations 
as we will show later by a truncation error analysis as well as by numerical experiments. 
Another way to see this accuracy deterioration for nonlinear equations is to recognize that the QUICKEST scheme is also equivalent to the conservative high-order finite-difference scheme of Shu and Osher in Ref.\cite{Shu_Osher_Efficient_ENO_II_JCP1989}. 
Based on the theory presented in Ref.\cite{Shu_Osher_Efficient_ENO_II_JCP1989}, the flux needs to be directly reconstructed instead of being evaluated with reconstructed solutions to achieve high-order accuracy. That is, we compute the numerical flux as
\begin{eqnarray}
F^{c}(u_L,u_R)  = \frac{1}{2} \left[  f_L   + f_R     \right]  - \frac{D}{2}  ( u_R - u_L), 
\label{upwnd_flux_convection_ENO}
\end{eqnarray}
where the left and right fluxes are directly reconstructed by 
\begin{eqnarray}
f_L &=& \frac{1}{2}  \left(  {f}_{i} +   {f}_{i+1}   \right)- \frac{1-\kappa}{4} \left(   {f}_{i+1}  -2 {f}_{i}+  {f}_{i-1} \right) ,  
\label{kappa_uL_simple_flux} \\ [2ex]
f_R &=& \frac{1}{2}  \left(  {f}_{i+1} +  {f}_{i}   \right)- \frac{1-\kappa}{4} \left(  {f}_{i+2}  -2  {f}_{i+1}+ {f}_{i} \right),
\label{kappa_uR_simple_flux}
\end{eqnarray}
with $\kappa=1/3$, where $f_i = f(u_i)$, and similarly for others. See Ref.\cite{Shu_WENO_SiamReview:2009} for further details. 
A more detailed discussion will be provided in a subsequent paper. 
Later, we will demonstrate numerically that third-order accuracy can actually be achieved for a nonlinear equation by the above direct flux reconstruction method.

To develop a genuinely third-order time-accurate QUICK scheme, we consider a simple technique to convert the cell-average time-derivative 
to a point-value time-derivative. In this approach, the residual stencil is kept within the five-point stencil, but it introduces a global coupling of the point-value time-derivatives, which needs to be
inverted at each explicit time step, as we will discuss in the next section.

\subsubsection{Coupled QUICK scheme}


Consider evaluating the second-derivative $ (u_{xx})$ with the central difference formula,
\begin{eqnarray}
 \overline{u}_{i }  
 =  u_i + \frac{1}{24} (u_{xx}) h^2  + O(h^4) 
 = u_i + \frac{1}{24} \left( \frac{  {u}_{i+1}  -2 {u}_{i}+  {u}_{i-1}  }{h^2} \right) h^2  + O(h^4).
 \label{ca_exact_expansion}
\end{eqnarray}
Then, 
\begin{eqnarray}
 \frac{d   \overline{u}_i }{dt} =   \frac{d   {u}_i }{dt}  +  \frac{1}{24} \left(
    \frac{d   {u}_{i+1} }{dt}  - 2   \frac{d  {u}_i }{dt} +   \frac{d   {u}_{i-1} }{dt} 
   \right)  + O(h^4),
    \label{ca_exact_expansion_dudt}
\end{eqnarray}
by which the scheme (\ref{fv_exact_form_discretized}) becomes
\begin{eqnarray}
  \frac{d   {u}_i }{dt}  +  \frac{1}{24} \left(
    \frac{d   {u}_{i+1} }{dt}  - 2   \frac{d  {u}_i }{dt} +   \frac{d   {u}_{i-1} }{dt} 
   \right)   + \frac{1}{h}  [ F_{i+1/2}  - F_{i-1/2}  ]   = \overline{s}_i,  
      \label{pointwise_scheme000}
\end{eqnarray}
or
\begin{eqnarray}
  \frac{1}{24} \left(   \frac{d   {u}_{i+1} }{dt}
  + 22    \frac{d   {u}_i }{dt} 
   +   \frac{d   {u}_{i-1} }{dt}    \right) + \frac{1}{h}  [ F_{i+1/2}  - F_{i-1/2}  ]   = \overline{s}_i.
   \label{pointwise_scheme}
\end{eqnarray}
This scheme is third-order accurate with $\kappa=1/2$. A similar coupled scheme can be found in Ref.\cite{Denaro_CMAME:2015}. 
The coupled scheme was also recognized by Van Leer, but he never programmed it \cite{BVL_private_communication_2020}. 


The time-derivative coupling is common to discretizations with point-valued solutions.
For example, the continuous Galerkin method leads to a similar coupled time-derivatives and the matrix
of coupling is often called the mass matrix \cite{Zienkiewicz_book}. The third-order node-centered edge-based discretization method
also requires a special coupling formula to achieve third-order accuracy on arbitrary triangular and tetrahedral
grids \cite{nishikawa_liu_source_quadrature:jcp2017}. Also, the residual-distribution method is often formulated 
with a mass matrix for time-dependent problems \cite{Mario_abgrall_ExplicitRD:JCP2010}.
Note also that the coupling exists only in space and therefore both explicit and implicit time integration schemes can be used to discretize the time derivatives in time. 

For an explicit time-stepping scheme, a tridiagonal system needs to be solved before updating the solution
to the next time level (or to the next stage in a multi-stage scheme such as the Runge-Kutta methods).
It can be seen more clearly in a global vector form:
\begin{eqnarray}
{\bf M} \frac{d {\bf U}}{dt} + {\bf Res}({\bf U}) = {\bf 0},
   \label{pointwise_scheme_vector}
\end{eqnarray}
where $ {\bf U}$ is a global vector of numerical solutions stored over $n$ cells, $ {\bf Res}({\bf U}) $ denotes the vector of residuals (i.e., the
spatial discretization including the forcing term) over $n$ cells, and the matrix ${\bf M}$ is an $n$$\times$$n$ tridiagonal matrix, which may be 
called a mass matrix: 
\begin{eqnarray} \label{myeq}
{\bf M} = 
 \frac{1}{24} \left[ \begin{matrix}
22 &        1  &           &            & 1 \\
1   &      22  &     1    &             &   \\
     & \ddots & \ddots & \ddots &   \\
  & & 1& 22 & 1 \\
 1  &  &    & 1 & 22
 \end{matrix} \right],
\end{eqnarray}
where we have assumed a periodic condition, i.e., $u_0 = u_n$ and $u_{n+1} = u_1$. 
Then, one can formally invert the mass matrix, 
\begin{eqnarray}
 \frac{d {\bf U}}{dt} + {\bf M}^{-1} {\bf Res}({\bf U}) = {\bf 0},
   \label{pointwise_scheme_vector_inverted}
\end{eqnarray}
and then apply an explicit time-stepping scheme to integrate in time. That is, the mass matrix needs to be inverted at each time step or at each stage of a multi-stage scheme.
For an implicit time-stepping scheme, the residual is coupled with the solution at the next time level, and 
thus we will have to solve the entire system of nonlinear residual equations in the form (\ref{pointwise_scheme_vector}) at each time step.
In this paper, the coupled mass-matrix scheme will be referred to as the coupled QUICK scheme.

\section{Truncation Error Analysis}
\label{QUICK_I_TE}
 
Following the previous paper \cite{Nishikawa_3rdMUSCL:2020}, we derive the modified equation \cite{Hirt_ModifiedEq:JCP1968} 
by expanding the scheme with a Taylor series of a numerical solution and identify the truncation error as the leading term in the difference between the modified equation and the target equation. The expansion can be performed relatively easily than in the previous paper \cite{Nishikawa_3rdMUSCL:2020} because the numerical solution is a point value this time. Let us assume that the numerical solution is smooth and expand it in a 
Taylor series around a cell center $i$,
\begin{eqnarray}
u(x) = u_i + (u_x) (x-x_i)  + \frac{1}{2} (u_{xx}) (x-x_i)^2  + \frac{1}{6} (\partial_{xxx} u) (x-x_i)^3   + \frac{1}{24} (\partial_{xxxx} u) (x-x_i)^4 + O(h^5).
\label{point_taylor}
\end{eqnarray}
 It is very important to note here that the derivatives such as $(u_x)$ and $(u_{xx})$ are also point values defined at $x=x_i$, which can be easily verified by differentiating $u(x) $ and evaluating the result at $x=x_i$. We do not use the subscript $i$ for the derivatives to distinguish them from the finite-difference approximations (\ref{fd_aprpox_at_ip1}). 
It is emphasized that we will need to perform the final conversion of a point-valued equation to a cell-averaged equation since the target equation is the integral form (\ref{integral_form_operator}), as we have done in Ref.\cite{Nishikawa_3rdMUSCL:2020}. 
  
Below, we will examine the accuracy of each part in the following order: the interpolated solutions, the fluxes, the diffusive flux and 
the diffusion scheme, the total flux balance term, and then two unsteady QUICK schemes. As in the previous paper \cite{Nishikawa_3rdMUSCL:2020}, we will show as much detail as possible to leave no room for misunderstanding.

\subsection{Accuracy of interpolated solution}
\label{third_order_inter_sol_acc}
 
Let us expand the left and right solutions interpolated at the right face as in Equations (\ref{kappa_uL_simple}) and (\ref{kappa_uR_simple}): 
\begin{eqnarray}
u_L  &=& u_i + \frac{1}{2} (u_x) h +  \frac{\kappa}{4} (u_{xx}) h^2 +  \frac{1}{12} (u_{xxx}) h^3  +  \frac{1}{96} (u_{xxxx}) h^4+ O(h^5), \label{uL_expanded} \\ [1.5ex]
u_R &=&  u_i + \frac{1}{2} (u_x) h +  \frac{\kappa}{4} (u_{xx}) h^2 +  \frac{1}{2} \left(  \frac{\kappa}{4} - \frac{1}{6} \right) (u_{xxx}) h^3 
 + \frac{1}{8} \left(  \frac{7 \kappa}{6} - 1 \right) (u_{xxxx}) h^4+ O(h^5),   \label{uR_expanded} 
\end{eqnarray}
whose average is given by
\begin{eqnarray}
\frac{u_L+u_R}{2}
=  u_i + \frac{1}{2} (u_x) h +  \frac{\kappa}{4} (u_{xx}) h^2 +  \frac{1}{8} \left(   \kappa  - \frac{1}{3} \right)  (u_{xxx}) h^3 + 
 \frac{1}{4} \left(  \frac{\kappa}{3} - \frac{1}{4} \right)  (u_{xxxx}) h^4 + O(h^5),  
\end{eqnarray}
which matches the exact value expanded as
\begin{eqnarray}
u_{i+1/2}^{exact} =  {u}_{i}  + \frac{1}{2} (u_x) h  + \frac{1}{8} (u_{xx}) h^2   + \frac{1}{48}  (\partial_{xxx} u) h^3  
+ \frac{1}{384}  (\partial_{xxxx} u) h^4 + O(h^5),
\label{exact_reconstruction_uL}
\end{eqnarray}
up to the cubic term if we take $\kappa=1/2$, 
\begin{eqnarray}
\frac{u_L+u_R}{2}
=  u_i + \frac{1}{2} (u_x) h +  \frac{1}{8} (u_{xx}) h^2 +  \frac{1}{48}   (u_{xxx}) h^3 + 
 \frac{1}{48}   (u_{xxxx}) h^4 + O(h^5).
\end{eqnarray}
Therefore, the QUICK interpolation scheme constructs a cubic function exactly on uniform grids, when averaged over the left
and right values.
\newline
\newline
\noindent {\bf Remark}: As mentioned in Ref.\cite{Nishikawa_3rdMUSCL:2020} the MUSCL scheme, the cubic exactness is a special property of a quadratically exact algorithm on a regular grid; the quadratic exactness is sufficient to design a third-order scheme, at least for convection. The resulting one-order higher truncation error is the reason that the truncation error order matches the discretization error (i.e., solution error) order on regular grids; the truncation error order is typically one-order lower on irregular grids for convection. See Refs.\cite{nishikawa_liu_source_quadrature:jcp2017,Boris_Jim_NIA2007-08,DiskinThomas:ANM2010,Katz_Sankaran_JCP:2011} and references therein for further details. 
\newline
\subsection{Accuracy of flux}
\label{third_order_inter_flux_acc}
 
Next, we consider the expansion of the averaged flux and prove that the cubic exactness holds for the averaged flux as well, which is a critical step in the truncation error analysis for a nonlinear equation \cite{Nishikawa_3rdMUSCL:2020}.
Consider the following expansions,
\begin{eqnarray} 
f(u_L)   &=&  f(u_i)   + f_u  du_L + \frac{1}{2} f_{uu} du_L^2  + \frac{1}{6} f_{uuu} du_L^3 + O(du^4),\\ [1.5ex]
f(u_R)  &=& f(u_i)   + f_u  du_R + \frac{1}{2} f_{uu} du_R^2  + \frac{1}{6} f_{uuu} du_R^3 + O(du^4),
\end{eqnarray}
where all the derivatives, $f_u$,  $ f_{uu}$, and $ f_{uuu}$, are point values at $x=x_i$, and $du_L = u_L - u_i $ and $du_R = u_R- u_i $.  
Taking the average and expanding further with Equations (\ref{uL_expanded}) and (\ref{uR_expanded}), we obtain
 \begin{eqnarray}
 \frac{f(u_L) + f(u_R)  }{2} &=&
  f(u_i)   + \frac{1}{2} ( \partial_{u} f)(u_x) h  
  +    \frac{1}{8}  \left[  2  \kappa  f_u  (u_{xx})  +  f_{uu} (u_x)^2     \right]  h^2   \nonumber \\
  &+&  \frac{1}{48} \left[     f_{uuu}(u_x)^3 +   6  \kappa  f_{uu} (\partial_{x} u)  (u_{xx})    +  2 ( 3 \kappa - 1)  f_u (\partial_{xxx} u)        \right] h^3    + O(h^4) ,
\end{eqnarray}
which matches the exact flux at the face, expanded as
\begin{eqnarray}
f_{j+1/2}^{exact} 
 &=&  f(u_i)   + f_x \left( \frac{h}{2} \right)  + \frac{1}{2} f_{xx} \left( \frac{h}{2} \right)^2
 +   \frac{1}{6}  f_{xxx}    \left( \frac{h}{2} \right)^3  +  O(h^4) \\ [2ex]
&=&  f(u_i)   + \frac{1}{2} f_x  h  + \frac{1}{8} f_{xx}  h^2  +   \frac{1}{48}  f_{xxx}   h^3    + O(h^4),
\end{eqnarray}
up to the cubic term if we take $\kappa=1/2$, 
 \begin{eqnarray}
 \frac{f(u_L) + f(u_R)  }{2} &=&
  f(u_i)   + \frac{1}{2} f_u u_x h  + \frac{1}{8} \left[  f_u  u_{xx}  + f_{uu} (u_x)^2   \right] h^2 \nonumber \\
  &+&  \frac{1}{48} \left[        f_{uuu} (u_x)^3 +   3 f_{uu}  u_x  u_{xx}   
  +  f_u  u_{xxx}
                        \right] h^3    + O(h^4) \\ [2ex]
 &=& 
 f(u_i)   + \frac{1}{2} f_x  h  + \frac{1}{8} f_{xx}  h^2  +   \frac{1}{48}  f_{xxx}   h^3    + O(h^4).
\end{eqnarray}
Similarly, the averaged flux at the left face is also exact for a cubic flux. Therefore, the QUICK interpolation scheme leads to
a cubically exact flux at a face, when averaged over the left and right values. 

\subsection{Order of dissipation term}
\label{third_order_dissipation_ac}

As proved for a cell-averaged solution in Ref.\cite{Nishikawa_3rdMUSCL:2020}, the dissipation term $D (u_R-u_L)$ is of $O(h^4)$ for any $\kappa$, when subtracted over the right and left faces and thus will produce an $O(h^3)$ truncation error. Therefore, any second-order error term will be generated by the averaged flux. The same is true for a point-valued solution; the results are exactly the same as presented in Ref.\cite{Nishikawa_3rdMUSCL:2020}. For this reason, we do not consider the dissipation term in the analysis below.

\subsection{Accuracy of diffusive flux and diffusion scheme}
\label{third_order_diffusion_ac}

For the diffusive flux, the damping term can generate an $O(h^2)$ error because its coefficient has a factor $1/h$,
\begin{eqnarray}
  \nu  (u_x)_{i+1/2} =   \nu \frac{ (u_x)_L  +(u_x)_R }{2} 
  + \frac{\nu \alpha}{2 h} (u_R -u_L) 
=   \nu \frac{  {u}_{i+1} -  {u}_{i}  }{h}  + \frac{\nu \alpha}{2 h} (u_R -u_L).
\end{eqnarray}
Therefore, we expand both the averaged diffusive flux and the damping term,  
\begin{eqnarray}
  \nu  (u_x)_{i+1/2} 
  =\nu u_x +    \frac{ 1}{2}  \nu u_{xx}  h  + \frac{1}{8} \left[  \frac{4}{3}  + \alpha ( \kappa  -1 )   \right] \nu u_{xxx} h^2  
+  \frac{1}{48} \left[  2 + 3   \alpha (\kappa-1)   \right] \nu u_{xxxx} h^3
  + O(h^4),
\end{eqnarray}
which matches the exact derivative expansion,
\begin{eqnarray}
  \nu  (u_x)^{exact}_{i+1/2} 
&=&
   \nu (u_x)  +    \nu  u_{xx}  \left( \frac{h}{2} \right)  + \frac{ 1 }{2}   \nu u_{xxx}   \left( \frac{h}{2} \right)^2 
    + \frac{  1}{6}   \nu u_{xxxx}  \left( \frac{h}{2} \right)^3 + O(h^4) \\[2ex]
&=& \nu (u_x) + \frac{1}{2} \nu u_{xx} h + \frac{ 1}{8} \nu u_{xxx} h^2 + \frac{ 1}{48}\nu u_{xxxx} h^3 + O(h^4),
\end{eqnarray}
up to the cubic term if we have
\begin{eqnarray} 
  \frac{4}{3}  + \alpha ( \kappa  -1 )   = 1, \quad 2 + 3   \alpha (\kappa-1)  = 1,
\end{eqnarray}
both of which are satisfied with 
\begin{eqnarray} 
  \alpha = \frac{1}{3 (1-\kappa)},
  \label{alpha_condition}
\end{eqnarray}
thus giving
\begin{eqnarray}
  \nu  (u_x)_{i+1/2}  = 
\nu (u_x) + \frac{1}{2} \nu u_{xx} h + \frac{ 1}{8} \nu u_{xxx} h^2 + \frac{ 1}{48}\nu u_{xxxx} h^3 + O(h^4).
\end{eqnarray}
Similarly, the cubic exactness holds at the left face as well. 
Note that the exactness holds for an arbitrary $\kappa$ as long as we choose $\alpha$ by Equation (\ref{alpha_condition}): e.g., 
for the QUICK scheme ($\kappa=1/2$), we have
\begin{eqnarray} 
  \alpha =  \frac{2}{3}.
  \label{alpha_condition_QUICK}
\end{eqnarray}
This diffusion scheme corresponds to the scheme derived in Ref.\cite{Nishikawa_FANG_AQ:Aviation2020} and also to one of the diffusion schemes considered by Leonard in Ref.\cite{Leonard_AMM1995}, which can be seen clearly by writing the diffusive flux $\nu  (u_x)_{i+1/2}$  with $(\kappa,\alpha) = (1/2,2/3)$ as 
\begin{eqnarray}
\left.  \nu  (u_x)_{i+1/2} \right|_{\alpha=2/3}
 &=&  \nu  \left[ \frac{  {u}_{i+1} -  {u}_{i}  }{h}  + \frac{1}{3 h} (u_R -u_L) \right]  \\ [2ex]
 &=&  \nu  \left[ \frac{  {u}_{i+1} -  {u}_{i}  }{h}  - \frac{   u_{i+2} -3 u_{i+1}  + 3 u_i  - u_{i-1}   }{24}  \right] ,
\end{eqnarray}
which is identical to Equation (18) of Ref.\cite{Leonard_AMM1995}. This scheme is also equivalent to the one 
used in a high-order finite-volume method in Ref.\cite{Denaro:JCP2011}. It is interesting to look at the resulting diffusion scheme by substituting
Equation (\ref{alpha_condition}) into Equation (\ref{fv_exact_form_discretized_diff_kappa_alpha}):
 \begin{eqnarray}
  \frac{1}{h}  [ F^d_{i+1/2}  - F^d_{i-1/2} ]    =  - \nu \frac{  - u_{i-2} + 28 u_{i-1} -54 u_i + 28 u_{i+1} - u_{i+2}  }{24 h^2},
  \label{diffusion_scheme_correct_third_order}
\end{eqnarray}
where $\kappa$ has been cancelled out.
That is, the diffusion scheme compatible with the QUICK scheme is unique. 
It is also fourth-order accurate since 
it is expanded as
 \begin{eqnarray}
  \frac{1}{h}  [ F^d_{i+1/2}  - F^d_{i-1/2} ]  =  - \nu  u_{xx} - \frac{\nu u_{xxxx}}{24} h^2 + O(h^4),
  \label{diffusion_scheme_correct_third_order_expanded}
\end{eqnarray}
and the first two terms correspond to the cell-averaged diffusion operator $-\overline{\nu u_{xx}}$ up to $O(h^4)$, which the finite-volume scheme 
is designed to approximate. See below.

\subsection{Accuracy of the total flux balance term}
\label{third_order_total_flux_ac}

Collecting results from the previous sections, we expand the flux balance term as 
\begin{eqnarray}
  \frac{1}{h}  [ F_{i+1/2}  - F_{i-1/2}  ] &=&  f_x +  \frac{1}{24} \left[     f_{uuu} (u_x)^3 +  6 \kappa  f_{uu} u_x  u_{xx}   +  2 \left(
 3 \kappa - 1
  \right) f_u  u_{xxx}
\right] h^2    \nonumber \\ [1ex]
& &- \nu u_{xx}  -   \frac{\nu}{24} u_{xxxx} h^2  + O(h^3).
\label{flux_balance_expanded_general}
\end{eqnarray}
As mentioned earlier and repeatedly pointed out by Leonard \cite{LeonardMokhtari:IJNMF1990,Leonard_AMM1995}, the spatial operator that the finite-volume discretization is trying to approximate is not the pointwise $f_x- \nu u_{xx} $ but the cell average $\overline{f_x}-(\overline{\nu u_{xx}})$. 
Consider the cell average of $f_x- \nu u_{xx} $: 
\begin{eqnarray}
 (\overline{f_x})_i  -(\overline{\nu u_{xx}})_i   =
  \frac{1}{h} \int_{x_i - h/2}^{x_i+h/2} ( f_{x}   - \nu  u_{xx}) \, dx  =  f_{x}   - \nu  u_{xx} + \frac{1}{24} \left(  f_{x}   - \nu  u_{xx}  \right)_{xx} h^2   + O(h^4),  
\end{eqnarray}
from which we find
\begin{eqnarray}
 f_{x}   - \nu  u_{xx}  = (\overline{f_x})_i  -(\overline{\nu u_{xx}})_i - \frac{1}{24} \left(  f_{x}   - \nu  u_{xx}  \right)_{xx} h^2   + O(h^4).
\end{eqnarray}
Substituting it into (\ref{flux_balance_expanded_general}), we obtain
\begin{eqnarray}
  \frac{1}{h}  [ F_{i+1/2}  - F_{i-1/2}  ]  = 
   (\overline{f_x})_i  -(\overline{\nu u_{xx}})_i 
    +  \frac{1}{24} \left[      f_{uuu} (u_x)^3 +  6 \kappa  f_{uu} u_x  u_{xx}   +  2 \left(
 3 \kappa - 1
  \right) f_u  u_{xxx}  -  f_{xxx}
\right] h^2   + O(h^3),
\label{flux_balance_expanded_general2}
\end{eqnarray}
which becomes by $f_{xxx} =  f_{uuu} (u_x)^3 +  3 f_{uu} u_x  u_{xx}   +   f_u  u_{xxx}$
\begin{eqnarray}
  \frac{1}{h}  [ F_{i+1/2}  - F_{i-1/2}  ] =
   (\overline{f_x})_i  -(\overline{\nu u_{xx}})_i 
    +  \frac{2 \kappa -1}{8} \left[     f_{uu} u_x  u_{xx}   +   f_u  u_{xxx}  
\right] h^2   + O(h^3).
\label{flux_balance_expanded_general3}
\end{eqnarray}
The second-order error vanishes for $\kappa=1/2$, thus resulting in 
\begin{eqnarray}
  \frac{1}{h}  [ F_{i+1/2}  - F_{i-1/2}  ] =  (\overline{f_x})_i  -(\overline{\nu u_{xx}})_i   + O(h^3).
  \label{expanded_flux_balance_final}
\end{eqnarray}
That is, the flux balance term is third-order accurate. 
The conversion to the cell-averaged operator is exactly what is missing in the analysis in Ref.\cite{JohnsonMackinnon:CANM1992} and has made them
conclude that the QUICK scheme is second-order accurate. 

The analysis is not complete yet. At this point, we have just shown that the spatial discretization is third-order accurate. 
It remains to examine the spatial discretization of the time derivative term for the two approaches: 
the coupled QUICK and QUICKEST schemes. 

\subsection{Accuracy of the coupled QUICK scheme}
\label{third_order_quick_acc}

For the coupled QUICK scheme (\ref{pointwise_scheme}), substituting Equation (\ref{expanded_flux_balance_final}) into the 
scheme (\ref{pointwise_scheme}), we obtain 
\begin{eqnarray}
  \frac{1}{24} \left(   \frac{d   {u}_{i+1} }{dt}
  + 22    \frac{d   {u}_i }{dt} 
   +   \frac{d   {u}_{i-1} }{dt}    \right)     +  (\overline{f_x})_i  -(\overline{\nu u_{xx}})_i   + O(h^3) = \overline{s}_i.
  \label{TE_INTEGRAL_FORM_coupled_00}
\end{eqnarray}
Since the time derivative terms can be expanded as in Equation (\ref{ca_exact_expansion_dudt}), it can be written as
 \begin{eqnarray}
  \frac{d   \overline{u}_i }{dt}  + O(h^4)    + (\overline{f_x})_i  -(\overline{\nu u_{xx}})_i    + O(h^3) = \overline{s}_i,
  \label{TE_INTEGRAL_FORM_coupled_01}
\end{eqnarray}
and thus
 \begin{eqnarray}
  \frac{d   \overline{u}_i }{dt}      + (\overline{f_x})_i  -(\overline{\nu u_{xx}})_i  + O(h^3) = \overline{s}_i.
  \label{TE_INTEGRAL_FORM_coupled}
\end{eqnarray}
This is the modified equation for the coupled QUICK scheme.    By comparing with 
the target integral form (\ref{integral_form_operator}), we find that the truncation error is $O(h^3)$. 
Therefore, the coupled QUICK scheme is third-order accurate for a general nonlinear conservation law.

\subsection{Accuracy of QUICKEST scheme}
\label{third_order_quickest}

As discussed in Section \ref{quickest_explained}, the QUICKEST scheme is equivalent to the $\kappa=1/3$ scheme. 
Then, it cannot be third-order accurate because third-order is achieved only for $\kappa=1/2$ as we have proved in Section \ref{third_order_total_flux_ac}.

However, it can be third-order accurate as a finite-difference scheme for a linear convection equation.
To see this, let $f = a u$, where $a$ is a global constant, and the expansion of the flux balance (\ref{flux_balance_expanded_general}) becomes, 
since $f_{uu}  = f_{uuu}=0$ and the diffusion is ignored ($\nu=0$), 
\begin{eqnarray}
  \frac{1}{h}  [ F^c_{i+1/2}  - F^c_{i-1/2}  ] =  f_u u_x +  \frac{ 3 \kappa - 1}{12}  f_u  u_{xxx} h^2     + O(h^3)
  =  f_x +  \frac{ 3 \kappa - 1}{12} ( f_{xxx}) h^2     + O(h^3).
  \label{quickest_te_01}
\end{eqnarray}
which further becomes for $\kappa=1/3$,
\begin{eqnarray}
  \frac{1}{h}  [ F^c_{i+1/2}  - F^c_{i-1/2}  ] =  f_x   + O(h^3),
\end{eqnarray}
showing that the scheme is a third-order finite-difference scheme for the differential form, not the integral form.
Leonard pointed it out by comparing schemes for the cell-averaged operator and the differential operator 
and showed that the $\kappa=1/3$ scheme is third-order accurate for a steady linear convection 
equation (with the exact diffusion term) \cite{Leonard_AMM1994}. No discussion was given for nonlinear equations.

Today, it is well known that this scheme is a conservative finite-difference scheme of Shu and Osher \cite{Shu_Osher_Efficient_ENO_II_JCP1989},
and as such, it is high-order only on one-dimensional uniform 
grids \cite{Shu_WENO_SiamReview:2009,Merryman:JSC2003} and a flux reconstruction is required for nonlinear equations
as mentioned in Section \ref{quickest_explained}. In fact, it can be proved quite straightforwardly since we already have all necessary expansions.
Replacing $u$ by $f$ in all the expansions in Section \ref{third_order_inter_sol_acc}, we obtain the expansion of the averaged flux at the right face, i.e., 
the average of Equations (\ref{kappa_uL_simple_flux}) and (\ref{kappa_uR_simple_flux}), as
\begin{eqnarray}
\frac{f_L+f_R}{2}
=  f(u_i) + \frac{1}{2} (f_x) h +  \frac{\kappa}{4} (f_{xx}) h^2 +  \frac{1}{8} \left(   \kappa  - \frac{1}{3} \right)  (f_{xxx}) h^3 + 
 \frac{1}{4} \left(  \frac{\kappa}{3} - \frac{1}{4} \right)  (f_{xxxx}) h^4 + O(h^5),
\end{eqnarray}
and a similar expansion at the left face (the same but with a negative sign given to the odd order terms). Then, we find 
\begin{eqnarray}
  \frac{1}{h}  [ F^c_{i+1/2}  - F^c_{i-1/2}  ] =  f_x   +  \frac{ 3 \kappa - 1}{12}  (f_{xxx}) h^2 + O(h^3),
\end{eqnarray}
for the flux reconstruction, which clearly shows that the scheme is a third-order approximation to the differential form $f_x$ with $\kappa=1/3$. 
Note that the above equation is valid for a general nonlinear conservation law whereas Equation (\ref{quickest_te_01}) is valid only for 
a linear conservation law. Later, we will numerically verify third-order accuracy of the flux reconstruction version of the QUICKEST scheme. 

Finally, it is emphasized that as the QUICKEST scheme is a finite-difference scheme, any forcing/source term, if present in a target equation, must be added to the QUICKEST scheme (\ref{pointwise_scheme1aaaa}) as a point-value at a cell center in order to preserve third-order accuracy:
\begin{eqnarray}
  \frac{d   {u}_i }{dt}  +  \frac{1}{h}  [  {F}^c_{i+1/2}  - {F}^c_{i-1/2}  ]   =  s(x_i).
      \label{pointwise_scheme1aaaa_forcing}
\end{eqnarray}
Also, if a diffusion term is present in a target equation, it must be discretized with a fourth-order finite-difference scheme, not a 
fourth-order finite-volume scheme, whose difference will be discussed in the next section.

\section{Confusing as Finite-Difference Scheme}
\label{QUICK_FD}

As we have seen, there is no confusion about third-order accuracy of the QUICK scheme as a finite-volume scheme with 
point-valued solutions. However, confusions will arise quickly if we re-interpret the QUICK scheme as a finite-difference scheme 
applied to the differential form (\ref{diff_form}). Below, we will discuss how confusing it is to construct a finite-difference
scheme based on the QUICK interpolation scheme.  

\subsection{Truncation error as a finite-difference scheme}
\label{QUICK_FD_TE}

If we consider a conservative finite-difference scheme,
\begin{eqnarray}
  \frac{d   {u}_i }{dt} + \frac{1}{h}  [ F_{i+1/2}  - F_{i-1/2}  ] = s_i,
  \label{wrong_fd_form}
\end{eqnarray}
where $s_i = s(x_i)$ and the numerical fluxes are computed by using the QUICK interpolation scheme ($\kappa=1/2$), then it is only second-order accurate. If we expand this scheme using Equation (\ref{flux_balance_expanded_general}) with $\kappa=1/2$, we obtain
\begin{eqnarray}
  \frac{d   {u}_i }{dt}  +  f_x     - \nu u_{xx} + \frac{1}{24}  \left(  f_{x}   - \nu  u_{xx}  \right)_{xx}       h^2     + O(h^3) = s(x_i).
  \label{wrong_fd_form01}
\end{eqnarray}
The second-order error remains and therefore the finite-difference scheme (\ref{wrong_fd_form}) is a second-order scheme 
for the differential form (\ref{diff_form}). Note that this scheme is not the QUICK scheme originally introduced in Ref.\cite{Leonard_QUICK_CMAME1979}. 
Therefore, calling the scheme (\ref{wrong_fd_form}) the QUICK scheme causes confusions.

To achieve third-order accuracy without changing the flux balance term, one must modify the discretizations of the time-derivative and forcing terms 
to generate second-order errors that cancel the above second-order error. In general, it is not very easy to devise such compatible discretizations. 
However, it is very simple if one follows the finite-volume formulation: one would notice immediately that both terms must be cell-averaged. 
For example, we replace $s_i$ by
\begin{eqnarray}
\overline{s}_i = s_i  + \frac{1}{24} \frac{  s_{i-1} -2 s_i + s_{i+1} }{h^2} h^2 + O(h^4),
\label{expanded_forcing}
\end{eqnarray}
and similarly 
\begin{eqnarray}
  \frac{d   \overline{u}_i }{dt} =   \frac{d   {u}_i }{dt} + \frac{1}{24} \frac{ \displaystyle \frac{d   {u}_{i+1} }{dt} -2 \frac{d   {u}_i }{dt}  + \frac{d   {u}_{i-1} }{dt}  }{h^2} h^2 + O(h^4), 
\label{expanded_dudt}
\end{eqnarray}
to construct the following finite-difference scheme:
\begin{eqnarray}
  \frac{d   {u}_i }{dt} + \frac{1}{24} \left( \frac{d   {u}_{i+1} }{dt} -2 \frac{d   {u}_i }{dt}  + \frac{d   {u}_{i-1} }{dt}  \right)+ \frac{1}{h}  [ F_{i+1/2}  - F_{i-1/2}  ] =
   s_i  + \frac{1}{24} \frac{  s_{i-1} -2 s_i + s_{i+1} }{h^2} h^2.
  \label{wrong_fd_form_third}
\end{eqnarray}
This scheme is third-order accurate. Expanding the scheme, we obtain
 \begin{eqnarray}
  \frac{d   {u}_i }{dt}    + \frac{1}{24}  \left( \frac{d u}{dt} \right)_{\!\! xx} h^2  + 
   f_x  - \nu u_{xx}+ \frac{1}{24} f_{xxx} h^2  = s_i  + \frac{1}{24}  s_{xx} h^2   + O(h^3),
\end{eqnarray}
which can be factored as
 \begin{eqnarray}
  \frac{d   {u}_i }{dt}      + 
   f_x   - \nu u_{xx} - s_i =   -     \frac{1}{24} \left(       \frac{d {u} }{dt}  + f_x  - \nu u_{xx} - s  \right)_{\!\! xx}  h^2 + O(h^3).
   \label{TE_FD_01}
\end{eqnarray} 
The left hand side is the differential form approximated by the finite-difference scheme (\ref{wrong_fd_form_third}).
It looks like that there is a second-order error on the right hand side, but it vanishes. There are two ways to see it. 
 
 For one, we see the above equation as a modified equation satisfied by a smooth numerical solution and differentiate it twice,
 \begin{eqnarray}
 \left(  \frac{d   {u}_i }{dt}      + 
   f_x  - \nu u_{xx} - s_i  \right)_{xx} =   -     \frac{1}{24} \left(       \frac{d {u} }{dt}  + f_x  - \nu u_{xx} - s  \right)_{\!\! xxxx}  h^2 + O(h^3),
\end{eqnarray} 
and substitute it into the right hand side of Equation (\ref{TE_FD_01}) to get 
\begin{eqnarray}
  \frac{d   {u}_i }{dt}      + 
   f_x  - \nu u_{xx} - s_i =  O(h^3),
   \label{TE_FD_02}
\end{eqnarray} 
showing that the finite-difference scheme (\ref{wrong_fd_form_third}) is a third-order approximation to the differential form. 
This error cancellation mechanism is well known and exploited in some discretization methods. For example, the residual-based compact method 
of Lerat and Corre \cite{Corre_Lerat_2005} achieves high-order accuracy by constructing the residual deliberately designed to achieve the error cancellation. 
Other examples include the residual-distribution method is another example \cite{Mario_abgrall_ExplicitRD:JCP2010,nishikawa_roe:ICCFD_2004} and 
the third-order edge-based method \cite{nishikawa_liu_source_quadrature:jcp2017}.

The other is to consider the smooth function used to expand the scheme as an exact solution, satisfying 
the target equation $  \frac{d {u} }{dt}  + f_x  - \nu u_{xx} - s  = 0$. Then, the $O(h^2)$ term in Equation (\ref{TE_FD_01}) vanishes because the second-derivative
of the target equation vanishes; it then leaves a third-order truncation error. 
In either case, we see that the finite-difference scheme (\ref{wrong_fd_form_third}) is a third-order finite-difference scheme for the differential form.
 
As one can see, the scheme (\ref{wrong_fd_form_third}) is nothing but the coupled QUICK scheme.
Calling the scheme (\ref{wrong_fd_form_third}) the QUICK finite-difference scheme will generate confusions. 
Only those having point-valued time-derivative and forcing terms, as in Equation (\ref{wrong_fd_form}), should be called the finite-difference scheme. 
As one might have noticed by now, it is very confusing to look at the QUICK scheme as a finite-difference scheme and 
 therefore not recommended. In the rest of the section, we will discuss some specific examples.

\subsection{Steady convection}

If we consider a steady conservation law $f_x=0$, then we obtain from Equation (\ref{TE_FD_01})
\begin{eqnarray}
   f_x  =   -  \frac{1}{24}   f_{xxx}  h^2 + O(h^3),
   \label{steady_fx_te}
\end{eqnarray}
which has been frequently used to argue that the QUICK scheme is second-order accurate. 
However, any discussion on accuracy is inconclusive here because $f_x=0$ implies that the exact solution
is given by a constant flux,
\begin{eqnarray}
   f  = constant,
\end{eqnarray}
and therefore any consistent discretization will be exact. All high-order derivatives in the error terms will be zero in Equation (\ref{steady_fx_te}). 
To verify the order of accuracy, one must 
include at least a forcing term $f_x=s(x)$. As shown in the previous section, the forcing term must be
cell-averaged in order to achieve third-order accuracy.

\subsection{Convection diffusion}
\label{subsec:confusion_advdiff}

Next, we consider Leonard's third-order numerical results shown in Ref.\cite{Leonard_AMM1995},
where he presents numerical results obtained for a one-dimensional steady convection-diffusion equation,
showing third-order accuracy with the QUICK scheme. 
He points out that the difference approximation to the diffusion term needs to be correctly chosen to achieve third-order accuracy with the QUICK convection 
scheme and presents the scheme that we rediscovered in Section \ref{third_order_diffusion_ac}: the alpha-damping scheme with 
$\alpha=2/3$ and $\kappa=1/2$.  Note that there is a typo in Equation (17) of Ref.\cite{Leonard_AMM1995}: 
the coefficient for $\phi_i$ should be $54$, not $52$. 
If viewed as a finite-difference scheme, this scheme has the following modified equation for a convection-diffusion 
equation:
 \begin{eqnarray}
   f_x   - \nu u_{xx} =   -     \frac{1}{24} \left(  f_x   - \nu u_{xx}   \right)_{xx}  h^2 + O(h^3) 
  \quad  \rightarrow  \quad 
   f_x  - \nu u_{xx}  =   O(h^3) .
   \label{compatible_diff_modified_eq_full}
\end{eqnarray}
As a finite-difference scheme, this particular diffusion scheme generates a second-order error, but it is canceled by the second-order error
coming from the QUICK scheme as mentioned earlier, thus resulting in a third-order convection-diffusion scheme.  

Hence, if the QUICK scheme is used as a finite-difference scheme, the diffusion scheme needs to be carefully chosen in order to achieve
third-order accuracy. This important point seems to have been often missed.  
In Ref.\cite{JohnsonMackinnon:CANM1992}, numerical results are presented for showing that the QUICK scheme is second-order 
accurate. The results were obtained for a linear convection-diffusion equation with the QUICK scheme combined with the
following fourth-order central difference diffusion scheme:
 \begin{eqnarray}
\frac{1}{h}  [ F^d_{i+1/2}  - F^d_{i-1/2} ]   = -  \nu \frac{  - u_{i-2} + 16 u_{i-1} -30 u_i + 16 u_{i+1} - u_{i+2}  }{12 h^2},  
  \label{diffusion_scheme_wrong}
\end{eqnarray}
which is expanded as
\begin{eqnarray}
\frac{1}{h}  [ F^d_{i+1/2}  - F^d_{i-1/2} ]   = -  \nu  u_{xx}  + \frac{u_{xxxxx}}{90} h^4 + O(h^5),
\end{eqnarray}
thus not generating any second-order error term and leading to the following modified equation for the full scheme,
 \begin{eqnarray}
   f_x   - \nu u_{xx} =   -     \frac{1}{24} \left(  f_x    \right)_{xx}  h^2 + O(h^3).
\end{eqnarray}
The second-order error remains and therefore the scheme is second-order. So, they did not use a compatible diffusion scheme. 
This diffusion scheme is suitable for a finite-difference scheme such as the QUICKEST scheme (\ref{pointwise_scheme1aaaa}), but not for a finite-volume scheme.
Leonard published a rebuttal to their claim in Ref.\cite{Leonard_AMM1995}. He recognized that the 
target diffusion operator was a cell-averaged operator, not a point-valued operator, and presented the diffusion scheme
equivalent to the one derived in this paper, without a derivation nor a truncation error analysis, which we have provided 
in Section \ref{third_order_inter_flux_acc}. As we have shown, this diffusion scheme generates a compatible second-order error
and successfully achieves third-order accuracy when combined with the QUICK scheme. Note that the scheme (\ref{diffusion_scheme_wrong})
can be derived from the alpha-damping scheme with $\kappa=1/2$ and $\alpha=4/3$, which does not satisfy the condition for high-order accuracy 
(\ref{alpha_condition}) and thus cannot be high-order.

The point here is that it is so much easier and straightforward to stick to the finite-volume formulation and construct a diffusion scheme
as a finite-volume scheme. Then, there is no confusion.

\begin{table}[t]
\ra{1.5}
\begin{center}
{\tabulinesep=1.2mm
\begin{tabu}{rrr}\hline\hline 
\multicolumn{1}{r}{ }                                               &
\multicolumn{1}{c}{ Finite-volume (integral form) }                                               &
\multicolumn{1}{c}{ Finite-difference (differential form)}    
\\ 
\multicolumn{1}{r}{ }                                               &
\multicolumn{1}{c}{ $ \displaystyle \frac{d  \overline{u}_i }{dt}  +\frac{1}{h}  [ F_{i+1/2}  - F_{i-1/2} ] =   \overline{s}_i$ }                                               &
\multicolumn{1}{c}{ $ \displaystyle   \frac{d   {u}_i }{dt}  + \frac{1}{h}  [ F_{i+1/2}  - F_{i-1/2} ]  =   {s}_i$  }    
  \\ \hline 
{\it steady scheme } &    &    \\
QUICK ($\kappa=1/2$) &   $\boldsymbol { \color{red} O(h^3)}$  &    $O(h^2)$    \\
 \hline 
{\it unsteady schemes } &    &    \\
Coupled QUICK ($\kappa=1/2$) &  $\boldsymbol { \color{red} O(h^3)}$   & \mbox{Lumped} - $O(h^2)$   \\
QUICKEST ($\kappa=1/3$) &       $O(h^2)$  &      $\boldsymbol { \color{red} O(h^3)}$  
  \\  \hline  \hline  
\end{tabu}
}
\caption{Spatial orders of accuracy of various QUICK schemes on uniform grids in one dimension. The numerical solution is a point value in all the cases. 
The third-order QUICKEST scheme requires a flux reconstruction for nonlinear equations. 
}
\label{Tab.orders_of_accuracy}
\end{center}
\end{table}

\section{Summary}
\label{third_order_summary}

Table \ref{Tab.orders_of_accuracy} summarizes the formal orders of accuracy for various QUICK schemes on uniform grids in one dimension.
All methods are based on the point-valued numerical solution. To avoid any confusion, the corresponding semi-discrete equations are shown 
for the finite-volume and finite-difference schemes. The steady scheme means the spatial discretization. 

As we have shown, the QUICK scheme is a third-order finite-volume discretization as a spatial discretization.
The coupled QUICK scheme achieves third-order accuracy for unsteady problems, again as a finite-volume scheme.
It is not applicable as a finite-difference scheme because it cannot be written in the finite-difference form (\ref{wrong_fd_form}): a single point-valued time derivative 
with a spatial discretization, $du_i /dt= \cdots$. Nevertheless, one may generate such a scheme by ignoring the coupling as
\begin{eqnarray}
 \frac{d   {u}_i }{dt}  \leftarrow \frac{1}{24} \left(   \frac{d   {u}_{i+1} }{dt}
  + 22    \frac{d   {u}_i }{dt} 
   +   \frac{d   {u}_{i-1} }{dt}    \right) ,
   \label{lumpd_mas_matrix}
\end{eqnarray}
which is often called the lumped mass matrix approach. The resulting scheme is a finite-difference scheme but only second-order accurate 
as shown in Section \ref{QUICK_FD_TE}. Later, second-order accuracy of the lumped mass matrix approach will be demonstrated numerically. 
Finally, the QUICKEST scheme is a third-order finite-difference scheme as 
we have discussed (with the flux reconstruction
for nonlinear equatons) and second-order accurate as a finite-volume scheme. 
 
It is strongly discouraged to think of these schemes in terms of cell-averaged solutions.
Doing so will only generate confusions. If we choose the cell-average as a numerical solution, then we are talking about 
the MUSCL scheme as discussed in Ref.\cite{Nishikawa_3rdMUSCL:2020}, not the QUICK scheme.

\section{Numerical Results}
\label{results}

\subsection{Steady convection problem with a forcing term}
\label{results_steady}

We begin with a steady problem for Burgers's equation in $x \in [0,1]$: 
\begin{eqnarray}
f_x =  s(x),
\end{eqnarray}
where $f = u^2/2$ with the forcing term, 
\begin{eqnarray}
s(x) =  2  \sin (2 x) \cos(  2 x ) ,
\end{eqnarray}
so that the exact solution is given by
\begin{eqnarray}
u(x) = \sin (2  x).
\end{eqnarray}
The forcing term is integrated exactly over each cell:
\begin{eqnarray}
  \overline{s}_i   =   \frac{1}{h} \int_{x_i-h/2}^{x_i+h/2} {s}(x) \, dx 
  =  \frac{1}{2h} \left[ 
    \cos^2(  h - 2 x_i  ) -    \cos^2(  h + 2 x_i  ) \right].
\end{eqnarray}
We solve the system of nonlinear finite-volume residual equations for the point-valued numerical solutions: 
${u}_3$,  ${u}_4$, $ \cdots$,  ${u}_{n-3}$, ${u}_{n-2}$, with the exact solution values are given and fixed at the left two cells $i=1$ and $2$, and at the 
right cells $i=n-1$ and $n$, in order to exclude boundary effects, which are beyond the scope of the present study.
An implicit solver based on the exact Jacobian of the first-order scheme is used to solve the residual equations.
See Ref.\cite{nishikawa_liu_jcp2018}, for example, for further details of the implicit solver for a one-dimensional finite-volume scheme. 
To verify the order of accuracy, we solve the steady problem with $\kappa=0$, $1/2$, and $1/3$, over a series of grids with with 15, 31, 63, 127 cells. The choice $\kappa=1/2$ corresponds to the QUICK scheme.

First, we verify the order of truncation error numerically by substituting the exact solution into the residual, which 
is defined at a cell center $i$ by
\begin{eqnarray}
Res_i =
   \frac{1}{h}  [ F_{i+1/2}  - F_{i-1/2}  ]   - \overline{s}_i,
\end{eqnarray}
and taking the $L_1$ norm over the cells. We consider both the point-valued exact solution and the cell-averaged exact solution and thus define the following
two truncation error norms:
\begin{eqnarray}
L_1( {\cal T}_p) = \frac{1}{n-4} \sum_{i=3}^{n-2} Res_i( \{ u_i^{exact} \} ), \quad
L_1( {\cal T}_c) = \frac{1}{n-4} \sum_{i=3}^{n-2} Res_i( \{ \overline{u}_i^{exact} \} ),
\end{eqnarray}
where $Res_i( \{ u_i^{exact} \} )$ is the residual with the point-valued exact solution substituted, 
and $Res_i( \{ \overline{u}_i^{exact} \} )$ is the residual with the cell-averaged exact solution substituted. 
The cell-averaged exact solution $ \overline{u}_i^{exact}$ is computed by exactly integrating the 
exact solution:
\begin{eqnarray}
 \overline{u}_i^{exact}   =   \frac{1}{h} \int_{x_i-h/2}^{x_i+h/2}  \sin (2  x) \, dx 
  =  \frac{1}{2h} \left[ 
    \cos (  h - 2 x_i  ) -    \cos (  h + 2 x_i  ) \right].
\end{eqnarray}
Note that $Res_i( \{ u_i^{exact} \} )$ is expanded, when $\kappa=1/2$,  as
\begin{eqnarray}
Res_i( \{ u_i^{exact} \} )
=
  (\overline{f_x})_i  -  \overline{s}_i  + O(h^3) = O(h^3) ,
\end{eqnarray}
where we have used $(\overline{f_x})_i  -  \overline{s}_i = 0$, which is true for an exact solution. 
Therefore, we expect $L_1( {\cal T}_p)$ to be third-order for $\kappa=1/2$. Similarly, from the analysis in Ref.\cite{Nishikawa_3rdMUSCL:2020},
$L_1( {\cal T}_c)$ is expected to be third-order for $\kappa=1/3$ since it is the third-order MUSCL scheme.

  \begin{figure}[t] 
    \centering
                \begin{subfigure}[t]{0.22\textwidth}
        \includegraphics[width=\textwidth]{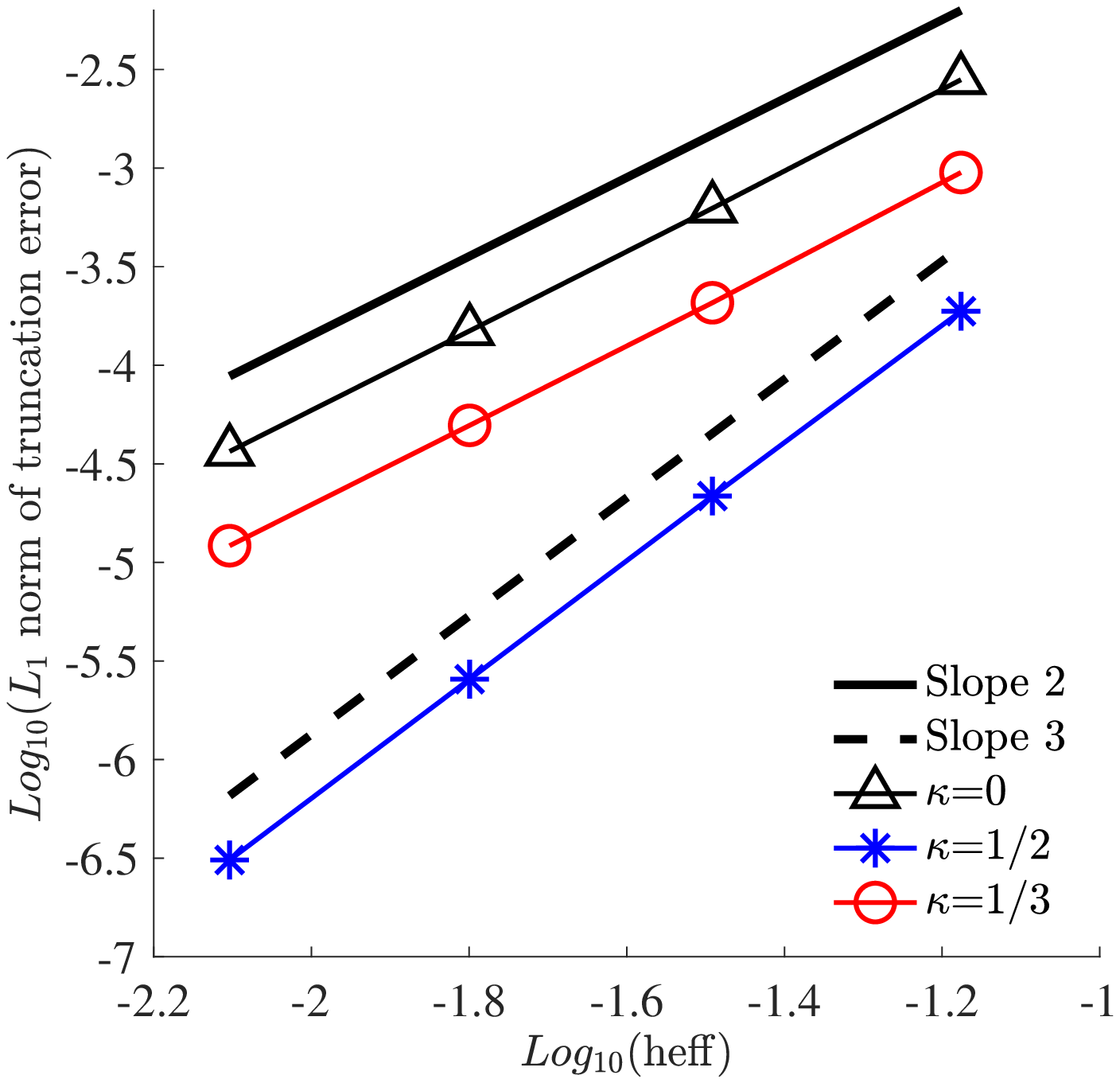}
          \caption{ $L_1( {\cal T}_p) $.}
       \label{fig:steady_adv_te_p}
      \end{subfigure}
      \hfill
          \begin{subfigure}[t]{0.22\textwidth}
        \includegraphics[width=\textwidth]{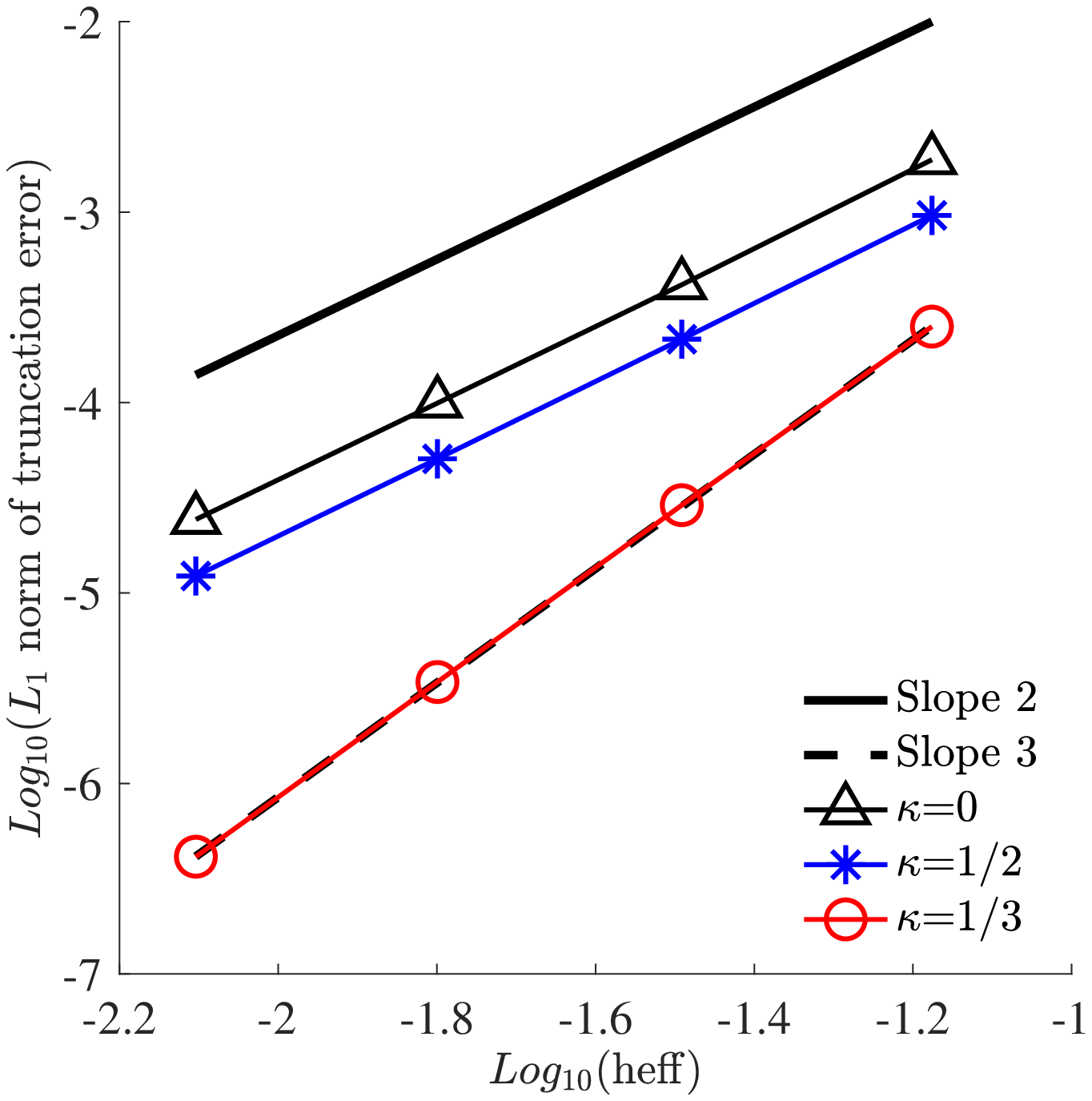}
          \caption{ $L_1( \hat{\cal T}_c ) $.}
       \label{fig:steady_adv_te_ca}
      \end{subfigure}
      \hfill
                \begin{subfigure}[t]{0.22\textwidth}
        \includegraphics[width=\textwidth]{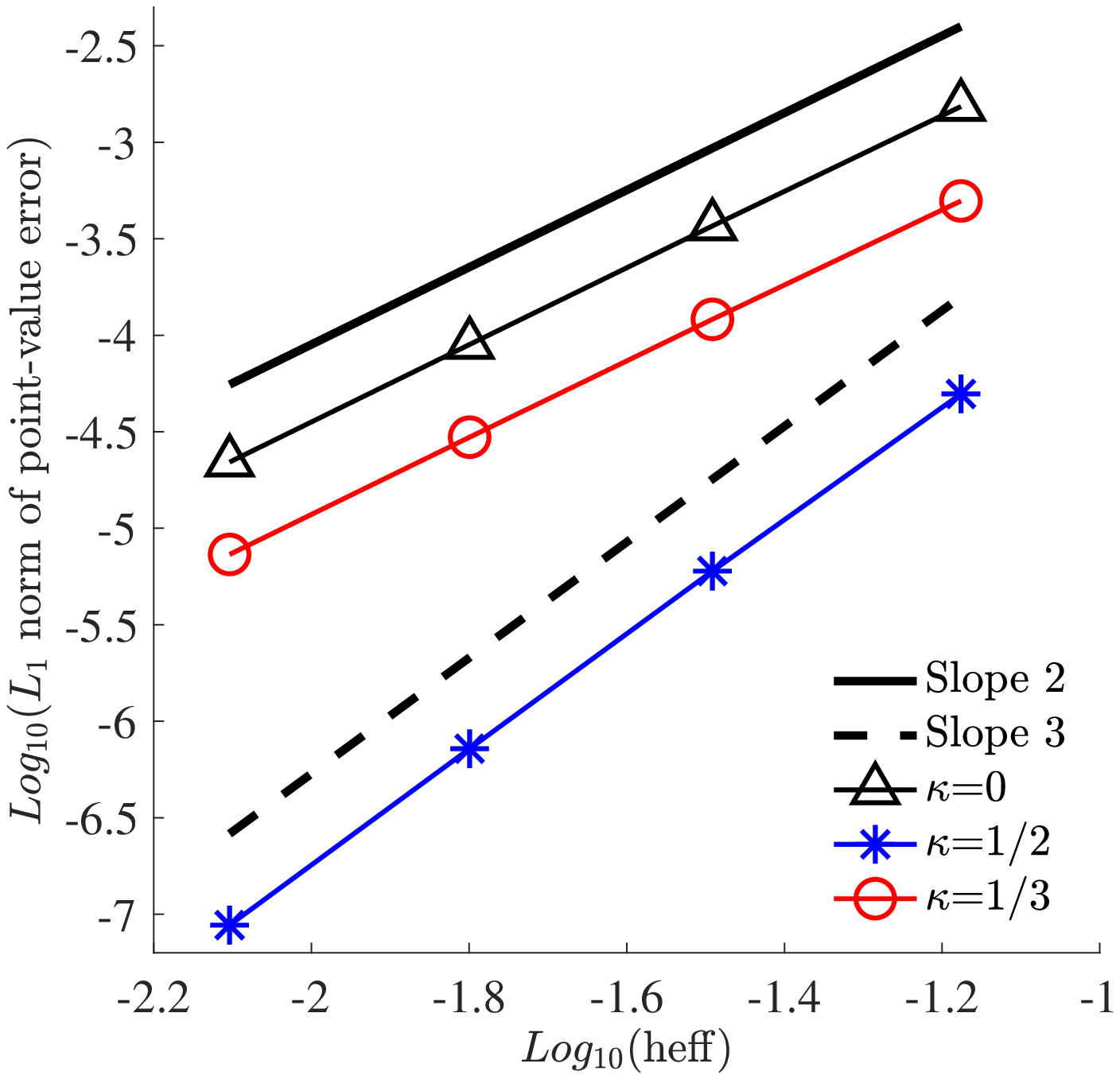}
          \caption{ $L_1( {\cal E}_p) $.}
       \label{fig:steady_adv_err_p}
      \end{subfigure}
      \hfill
          \begin{subfigure}[t]{0.22\textwidth}
        \includegraphics[width=\textwidth]{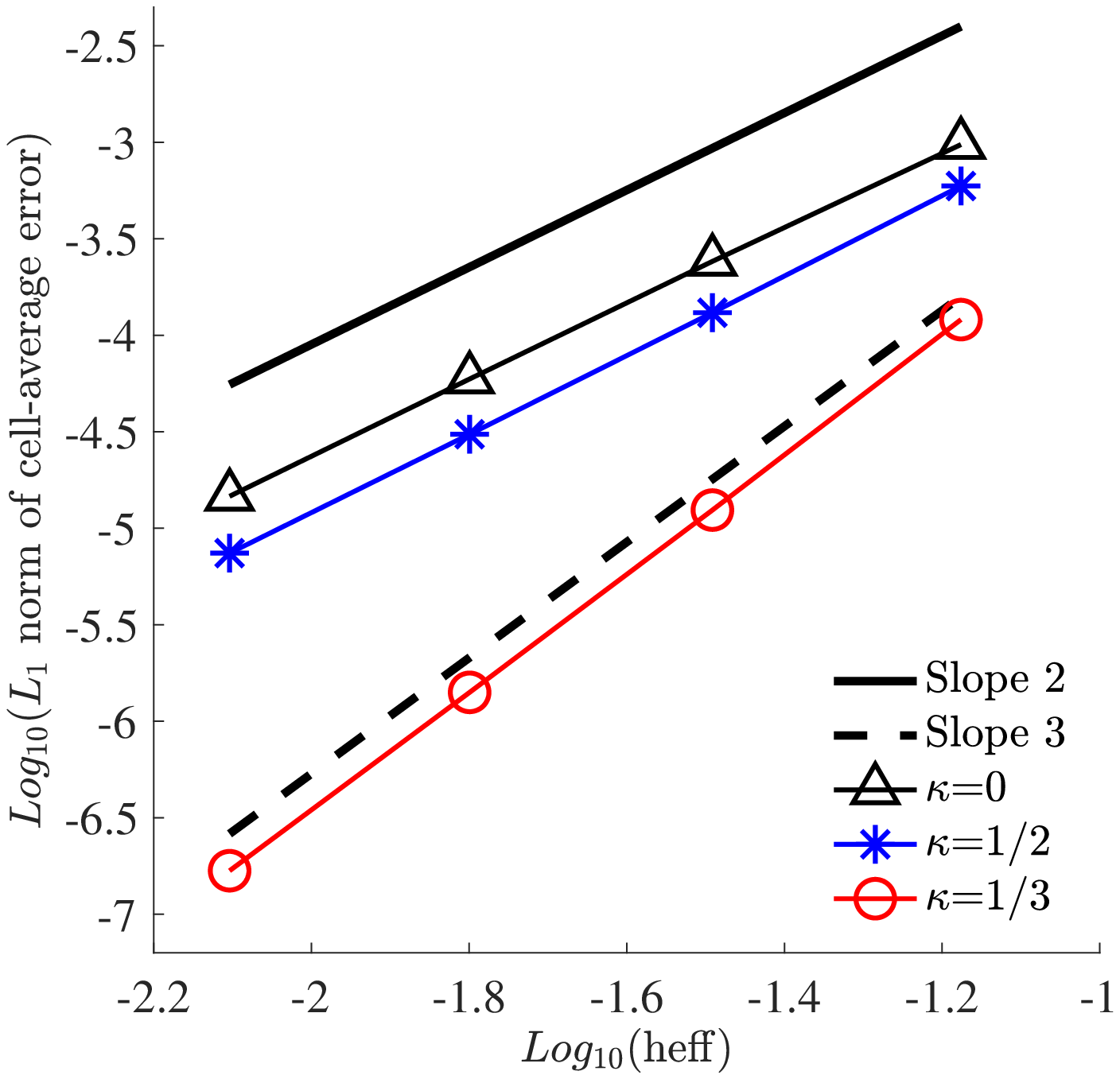}
          \caption{ $L_1( \hat{\cal E}_c ) $.}
       \label{fig:steady_adv_err_ca}
      \end{subfigure}
            \caption{
\label{fig:steady_adv_error}%
Truncation and discretization error convergence results for the steady Burgers equation.
} 
\end{figure}

Figure \ref{fig:steady_adv_te_p} shows the convergence of the truncation error $L_1( {\cal T}_p)$. 
 As expected, the truncation error $L_1( {\cal T}_p)$ is third-order for the QUICK scheme ($\kappa=1/2$). 
Figure \ref{fig:steady_adv_te_ca} shows the convergence of the truncation error computed with the cell-averaged exact solution, $L_1( {\cal T}_c)$.
As can be seen, third-order convergence is obtained with $\kappa=1/3$ as expected. These results verify that the steady residual equation can be 
considered as both the third-order MUSCL scheme ($\kappa=1/3)$ and the third-order QUICK scheme ($\kappa=1/2$), depending on how the numerical 
solution is interpreted, as discussed in Ref.\cite{Nishikawa_3rdMUSCL:2020}.

For the discretization error (i.e., solution error), we consider two error norms: the point-value error and the cell-average error,
\begin{eqnarray}
L_1( {\cal E}_p) = \frac{1}{n-4} \sum_{i=3}^{n-2} |  u_i - u_i^{exact} |,  \quad
L_1( {\cal E}_c) = \frac{1}{n-4} \sum_{i=3}^{n-2} | u_i - \overline{u}_i^{exact}|.
\end{eqnarray}
 Note that the numerical solution is considered as a point-valued solution in the QUICK scheme 
and thus always denoted by $u_i$.
Results are shown in Figures \ref{fig:steady_adv_err_p} and \ref{fig:steady_adv_err_ca}.
As expected, the QUICK scheme ($\kappa=1/2$) is third-order for the point-valued solution,
and the MUSCL scheme ($\kappa=1/3$) is third-order for the cell-averaged solution.

It is worth pointing out that the finite-volume scheme tested here is equivalent to the U-MUSCL scheme of
Burg \cite{burg_umuscl:AIAA2005-4999}. As stated in Ref.\cite{burg_umuscl:AIAA2005-4999}, Burg discretized a steady one-dimensional
shallow-water system by the node-centered finite-volume method with forcing terms integrated over a control volume. 
He then found that third-order accuracy was achieved for point-valued solutions at nodes with $\kappa=1/2$, which is exactly the QUICK scheme.
For this reason, his scheme should have been called the unstructured-QUICK (U-QUICK) scheme. Apparently, he did not recognize the difference between
point-valued and cell-averaged solutions and thus only provided a very brief heuristic argument about why third-order accuracy was not obtained with $\kappa=1/3$. 
Here, we have shown that the reason is the point-valued numerical solution, which makes the finite-volume scheme the QUICK scheme, not MUSCL. 
 
Finally, it is noted that third-order with $\kappa=1/3$ observed here is a special case for a pure convection equation 
and does not automatically carry over to a convection-diffusion problem, as we will
discuss in the next section.

\subsection{Steady convection-diffusion problem}
\label{results_steady_advdiff}
 
Next, we consider a steady problem for the viscous Burgers equation with a forcing term: 
\begin{eqnarray}
f_x =  \nu u_{xx} +s(x),
\end{eqnarray}
in $x \in [0,1]$, where $f = u^2/2$ and
\begin{eqnarray}
s(x) =  2  \sin (2 x) \cos(  2 x ) + 4 \nu \sin( 2 x),
\end{eqnarray}
so that the exact solution is given by
\begin{eqnarray}
u(x) = \sin (2  x).
\end{eqnarray}
Again, we integrate the forcing term is integrated exactly over each cell:
\begin{eqnarray}
  \overline{s}_i   =   \frac{1}{h} \int_{x_i-h/2}^{x_i+h/2} {s}(x) \, dx 
  =  \frac{1}{2h} \left[ 
    \cos^2(  h - 2 x_i  ) -    \cos^2(  h + 2 x_i  ) 
  \right] 
-  \frac{2 \nu }{h} \left[ 
    \cos(  h - 2 x_i  ) -    \cos(  h + 2 x_i  ) 
  \right].
\end{eqnarray}
In particular, we consider the case $\nu=1$, for which the convective and diffusion terms are equally important.
As before, we solve the steady problem by the implicit solver with $\kappa=0$, $1/2$, and $1/3$, 
over a series of uniform grids with 15, 31, 63, 127 cells. The damping coefficient $\alpha$ is determined by the
condition (\ref{alpha_condition}) and the same $\kappa$ is used as in the convective flux. Then, the analysis predicts that third-order accuracy 
is obtained only for $\kappa=1/2$.

  \begin{figure}[th!]
    \centering
                \begin{subfigure}[t]{0.3\textwidth}
        \includegraphics[width=\textwidth]{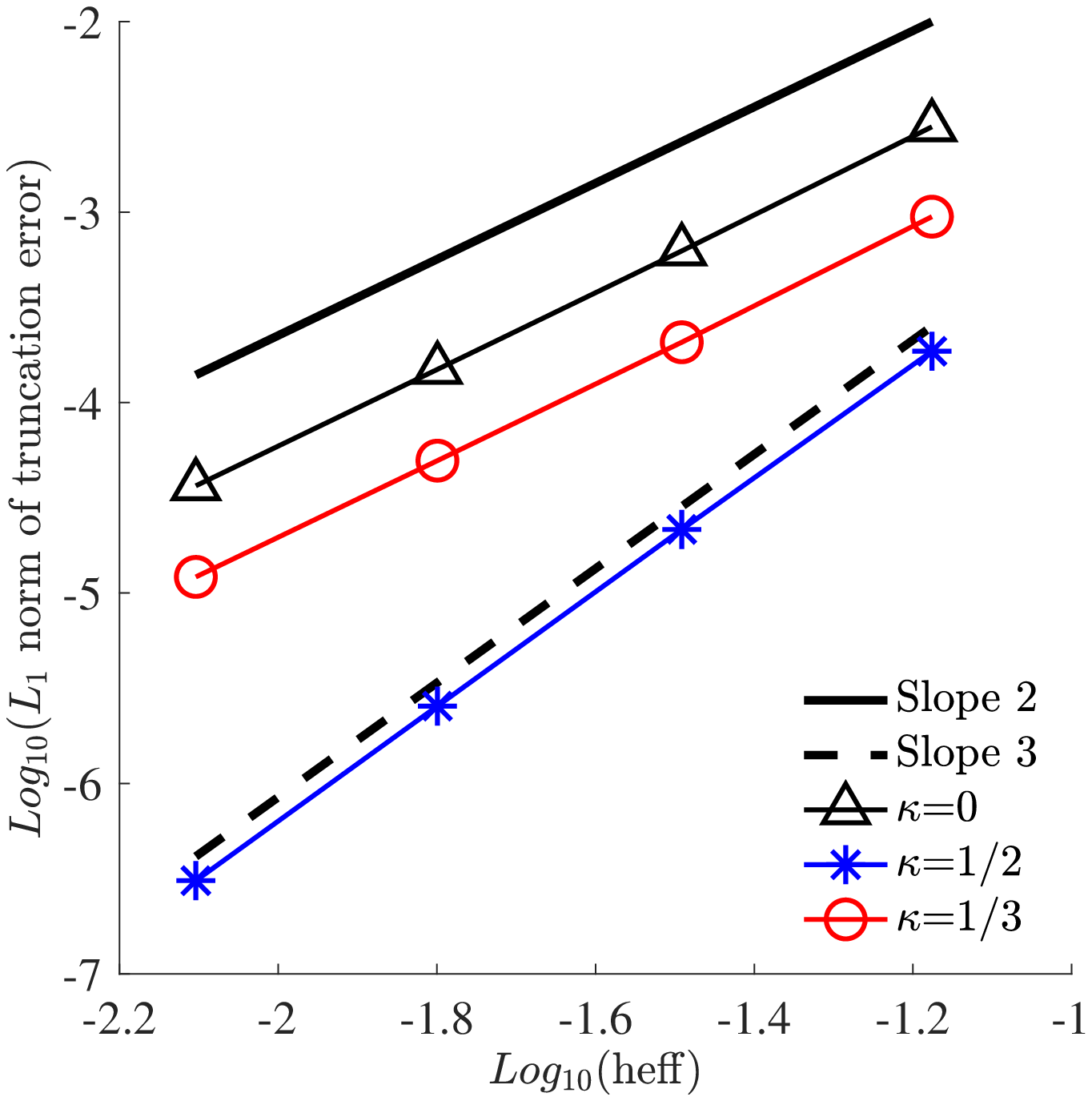}
          \caption{ $L_1( {\cal T}_p) $.}
       \label{fig:steady_adv_visc_te_p}
      \end{subfigure}
            \hfill
                \begin{subfigure}[t]{0.3\textwidth}
        \includegraphics[width=\textwidth]{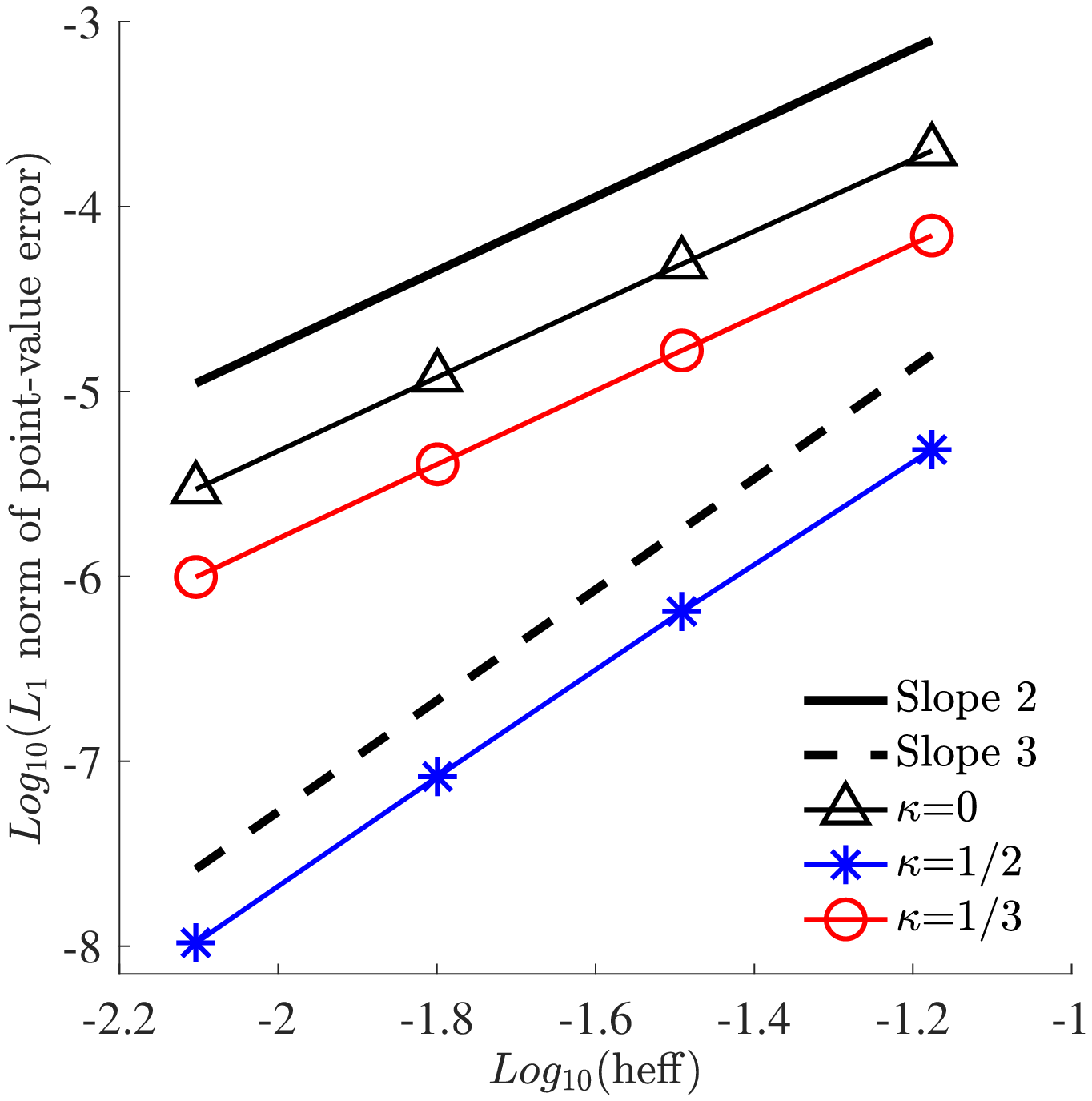}
          \caption{ $L_1( {\cal E}_p) $.}
       \label{fig:steady_adv_visc_err_p}
      \end{subfigure}
      \hfill
          \begin{subfigure}[t]{0.3\textwidth}
        \includegraphics[width=\textwidth]{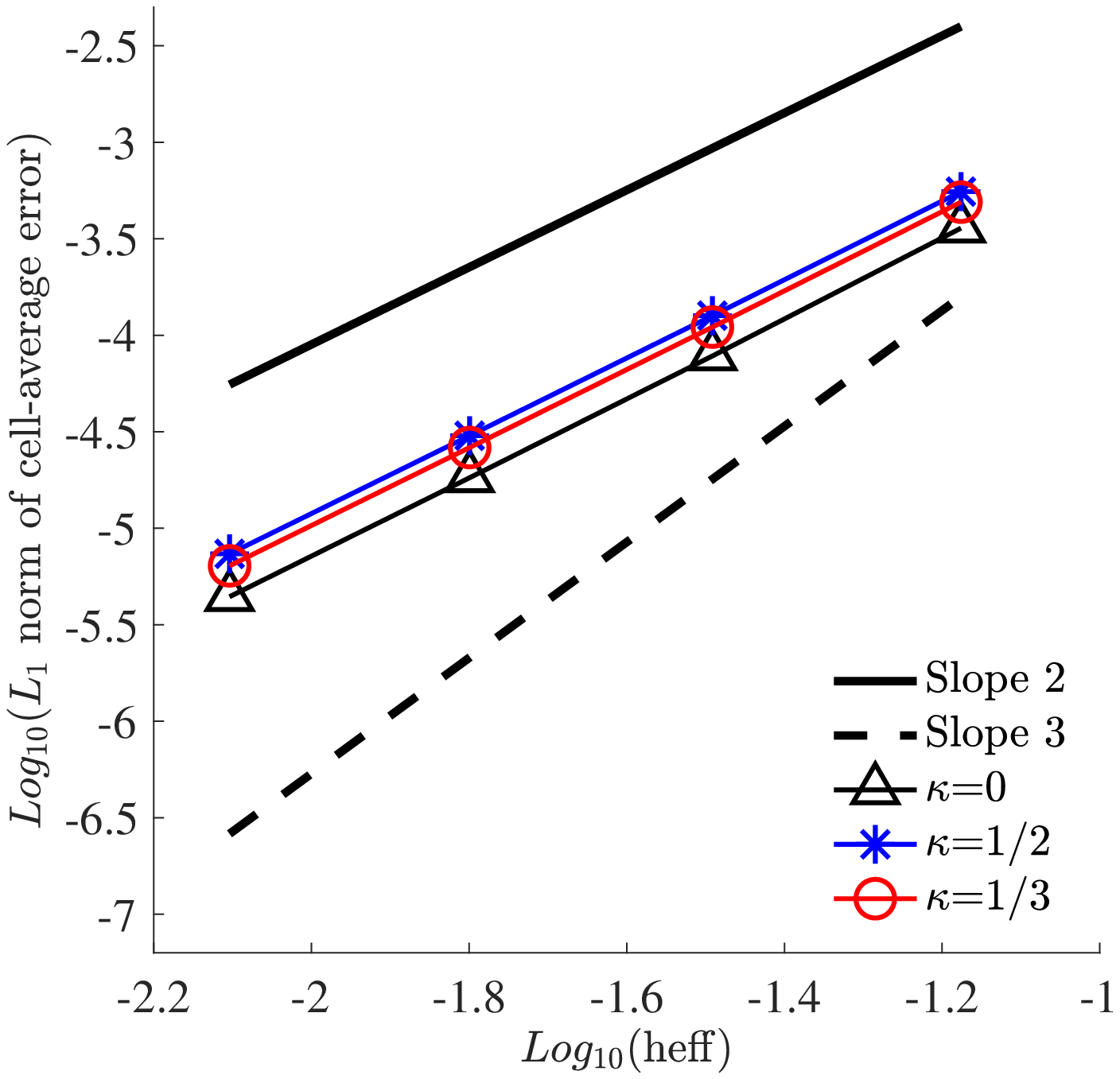}
          \caption{ $L_1( \hat{\cal E}_c ) $.}
       \label{fig:steady_adv_visc_err_ca}
      \end{subfigure}
            \caption{
\label{fig:steady_adv_visc_error}%
Truncation and discretization error convergence results for the case of the steady viscous Burgers equation.
The diffusion scheme is implemented as the alpha-damping scheme with $\alpha=\frac{1}{3 (1-\kappa)}$.
} 
\end{figure}
  \begin{figure}[th!]
    \centering
                \begin{subfigure}[t]{0.3\textwidth}
        \includegraphics[width=\textwidth]{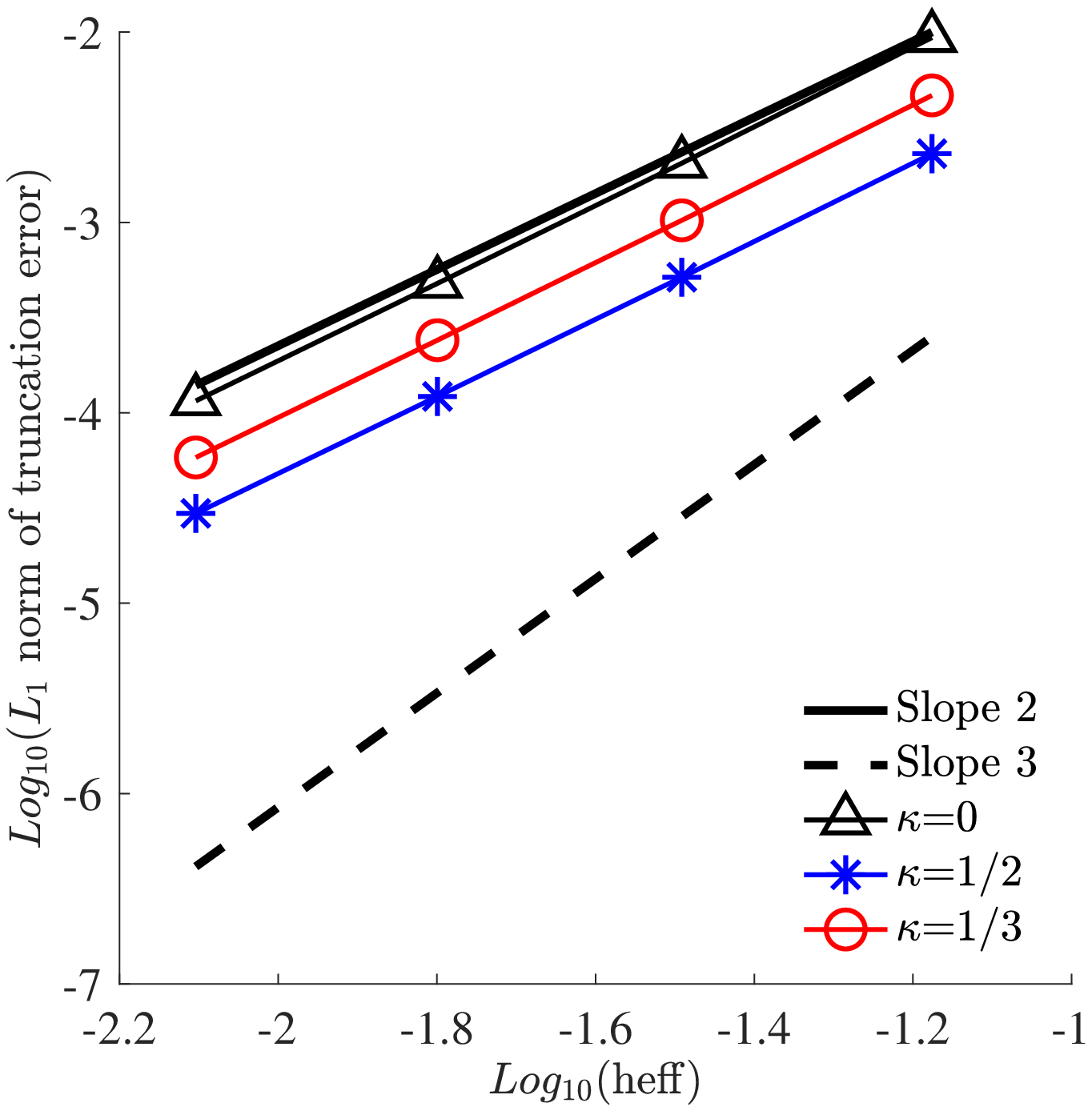}
          \caption{ $L_1( {\cal T}_p) $.}
       \label{fig:steady_adv_visc_te_p_alpha4o3}
      \end{subfigure}
            \hfill
                \begin{subfigure}[t]{0.3\textwidth}
        \includegraphics[width=\textwidth]{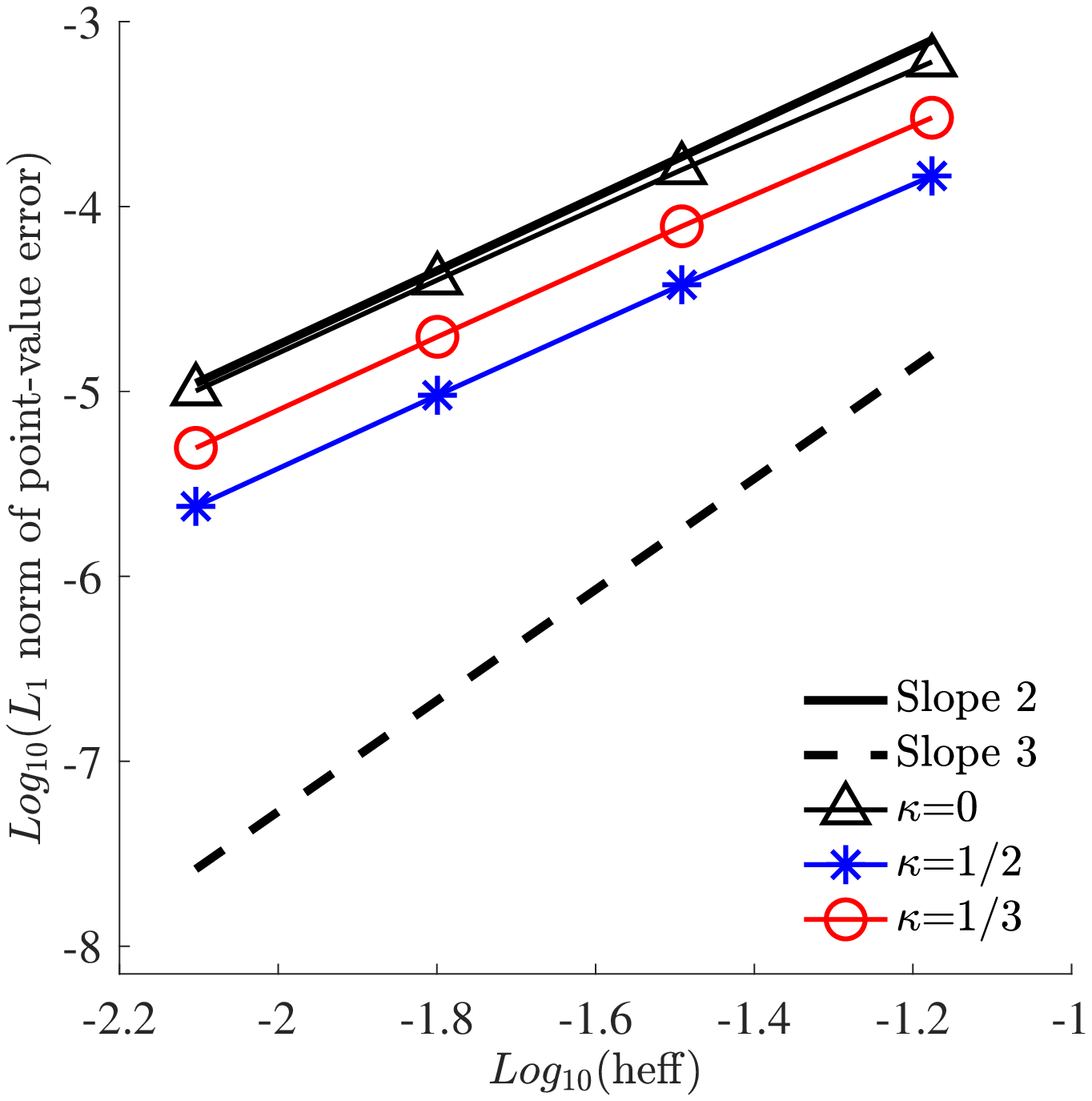}
          \caption{ $L_1( {\cal E}_p) $.}
       \label{fig:steady_adv_visc_err_p_alpha4o3}
      \end{subfigure}
      \hfill
          \begin{subfigure}[t]{0.3\textwidth}
        \includegraphics[width=\textwidth]{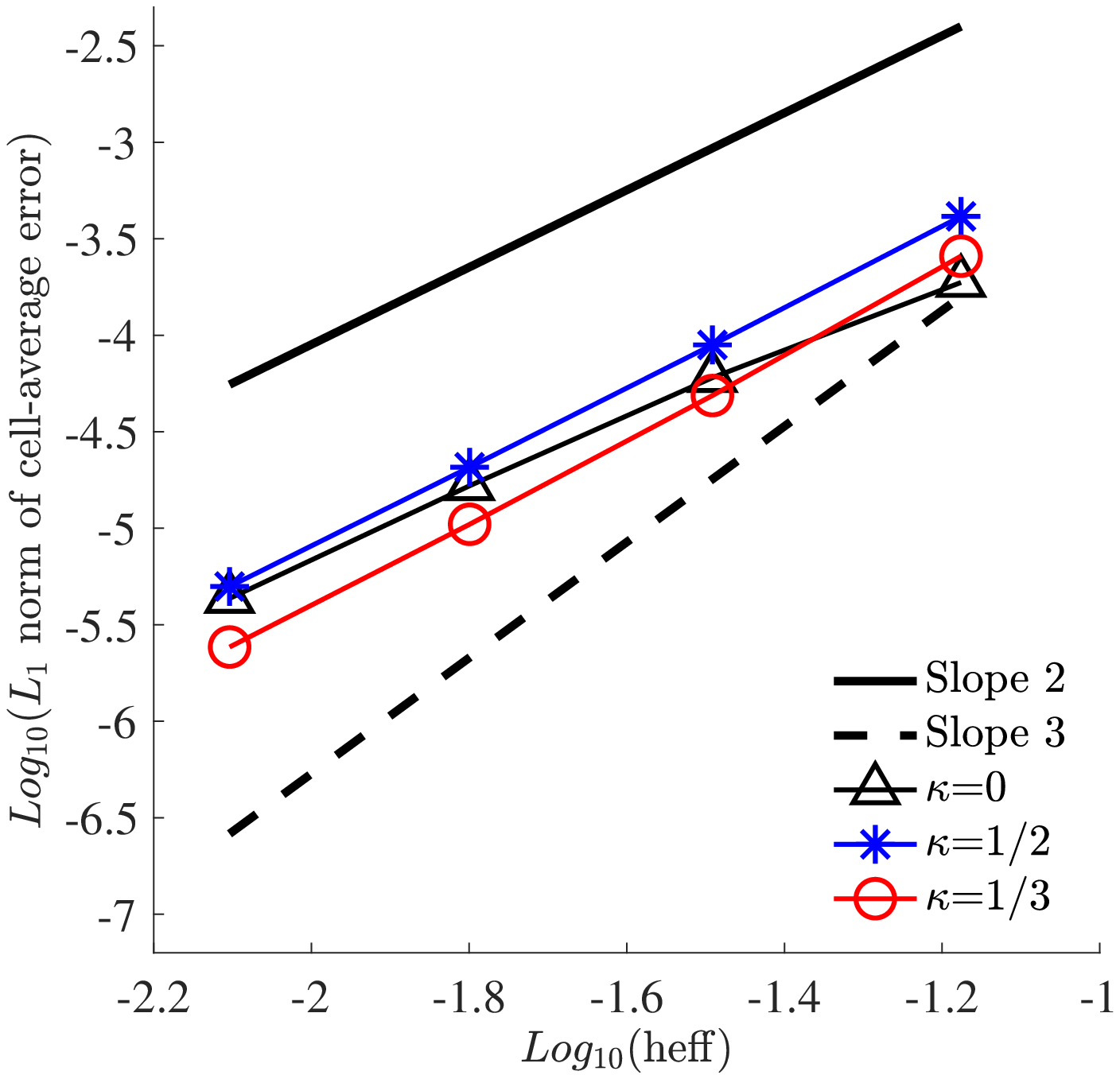}
          \caption{ $L_1( \hat{\cal E}_c ) $.}
       \label{fig:steady_adv_visc_err_ca_alpha4o3}
      \end{subfigure}
            \caption{
\label{fig:steady_adv_visc_error_alpha4o3}%
Truncation and discretization error convergence results for the case of the steady viscous Burgers equation. 
The diffusion scheme is implemented as the alpha-damping scheme with $\alpha=4/3$.
} 
\end{figure}
 
From this problem, we will focus on the QUICK scheme with point-valued solutions and consider only the truncation error $L_1( {\cal T}_p)$ 
and the discretization error $L_1( {\cal E}_p)$, but the error as the cell-averaged solution $L_1( {\cal E}_c)$ is also shown for clarity.
Figure \ref{fig:steady_adv_visc_te_p} shows the convergence of the truncation error $L_1( {\cal T}_p)$. 
As expected, third-order accuracy is achieved only with $\kappa=1/2$. Then, as we would expect, the discretization error $L_1( {\cal E}_p)$ is also third-order only with $\kappa=1/2$ as shown in Figure \ref{fig:steady_adv_visc_err_p}.
To emphasize that third-order accuracy is achieved for the point-valued solution, not for the cell-averaged solution, 
we present the discretization error convergence with the cell-averaged exact solution, $L_1( {\cal E}_c)$. See Figure \ref{fig:steady_adv_visc_err_ca}.
Clearly, third-order accuracy is not achieved for any value of $\kappa$. Thus, unlike the convection case in the previous section, third-order accuracy 
is not achieved automatically with $\kappa=1/3$ because the diffusion scheme is not fourth-order accurate with cell-averaged solutions. 
These results demonstrate again that the QUICK scheme is third-order accurate for the point-valued solution.

Finally, we consider the diffusion scheme (\ref{diffusion_scheme_wrong}), which was shown to yield second-order accuracy when
combined with the QUICK scheme in Section \ref{subsec:confusion_advdiff}. We implemented this scheme as the alpha-damping scheme 
with $\alpha=4/3$; then, the scheme (\ref{diffusion_scheme_wrong}) is reproduced when $\kappa=1/2$. 
Results are shown in Figure \ref{fig:steady_adv_visc_error_alpha4o3}. As expected, third-order accuracy is not observed 
for both the truncation and discretization errors.

\subsection{Unsteady convection problem} 
\label{subsec:unsteady}

\begin{figure}[t]
    \centering
      \hfill    
          \begin{subfigure}[t]{0.45\textwidth}
        \includegraphics[width=\textwidth]{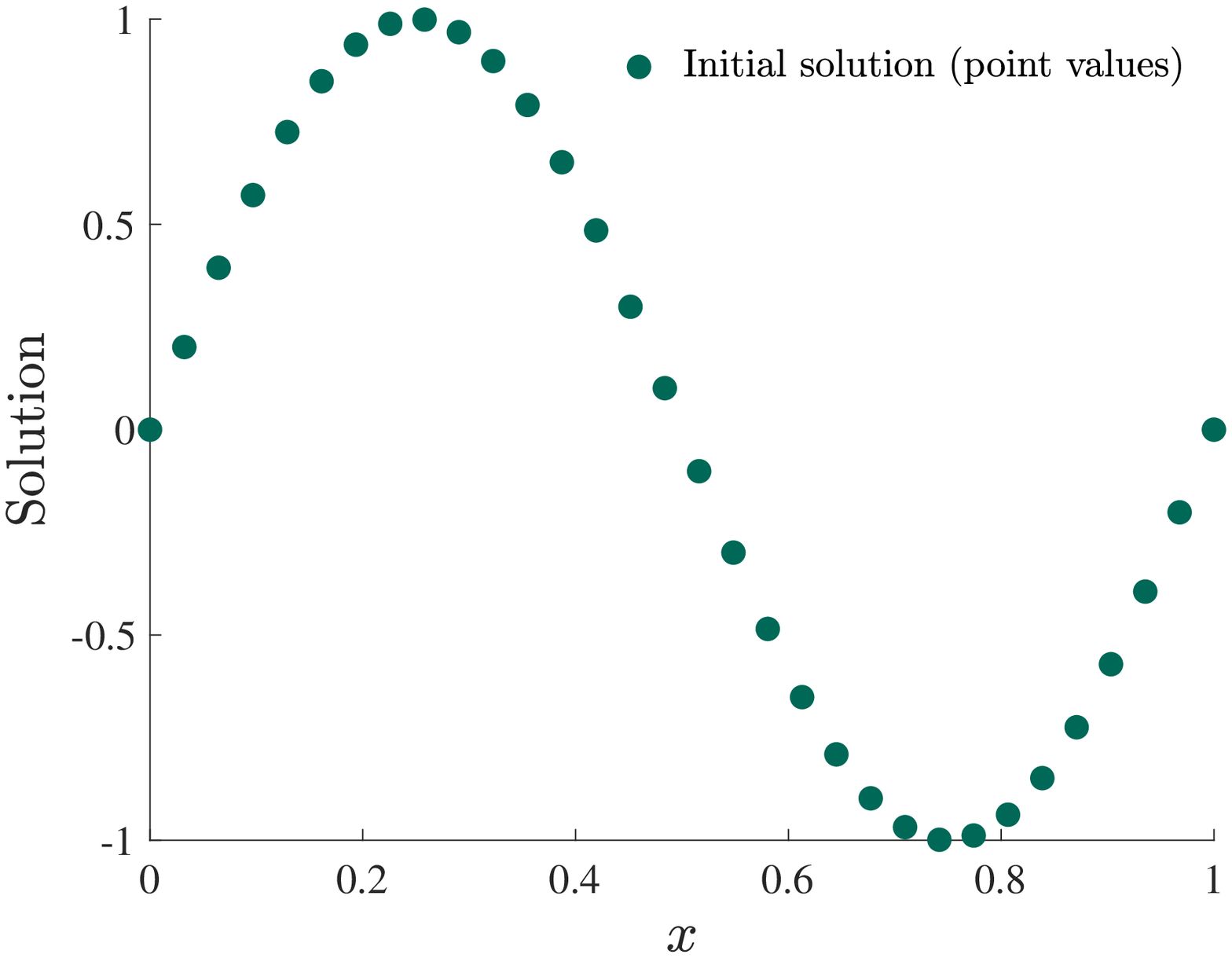}
          \caption{Initial solution (point values).}
          \label{fig:unsteady_initial}
      \end{subfigure}
      \hfill
          \begin{subfigure}[t]{0.45\textwidth}
        \includegraphics[width=\textwidth]{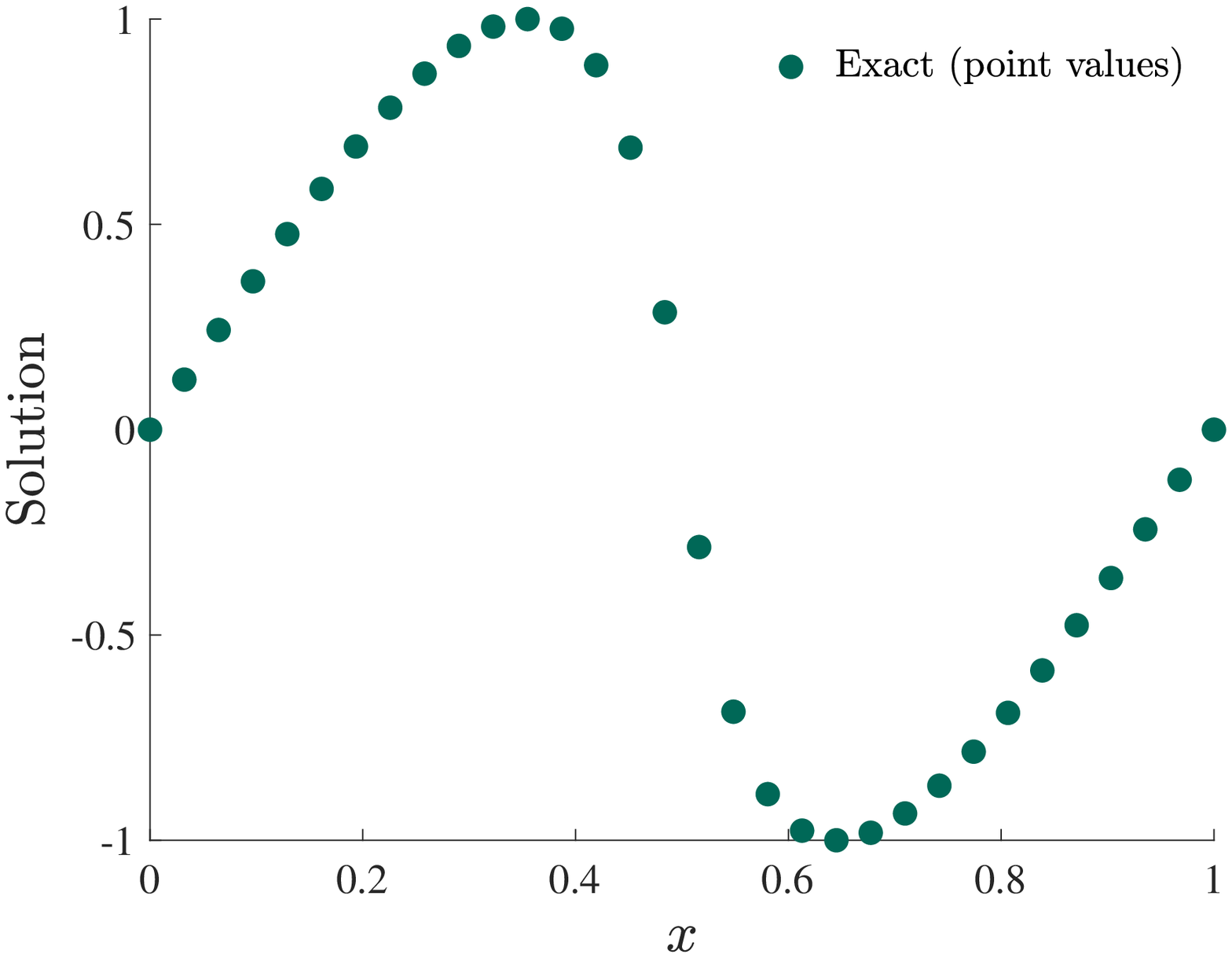}
          \caption{Final solution (point values).}
       \label{fig:unsteady_final}
      \end{subfigure}
      \hfill    
            \caption{
\label{fig:unsteady_solution}%
Initial and final solutions on the coarsest grid for the unsteady test case.
} 
\end{figure}

To demonstrate third-order accuracy for unsteady problems, we consider the same time-dependent problem for Burgers's equation as in the previous paper \cite{Nishikawa_3rdMUSCL:2020}: 
\begin{eqnarray}
u_t  + f_x= 0, 
\end{eqnarray}
where $f = u^2/2$, with the initial solution,
\begin{eqnarray}
u(x) = \sin (2 \pi x).
\end{eqnarray}
The domain is $x \in [0,1]$ but it is taken to be periodic (i.e., there is no boundary in this problem).
We compute the solution at the final time $t = t_f =0.105$. The initial and final solutions are shown in Figure \ref{fig:unsteady_solution}. 
Various schemes are compared: $\kappa=0$, $1/2$, and $1/3$ over a series of
grids: 32, 64, 128, 256, 512, 1024, 2048 cells. We implemented both the coupled QUICK scheme (with the mass matrix inverted directly by 
Thomas' algorithm; see, e.g., Ref.\cite{Hirsch_VOL1}) and the QUICKEST scheme, both of which are integrated in time by the three-stage SSP Runge-Kutta scheme \cite{SSP:SIAMReview2001} for the total of 840 time
steps with a constant time step $\Delta t = 0.000125$.
 
  \begin{figure}[t]
    \centering
          \begin{subfigure}[t]{0.3\textwidth}
        \includegraphics[width=\textwidth]{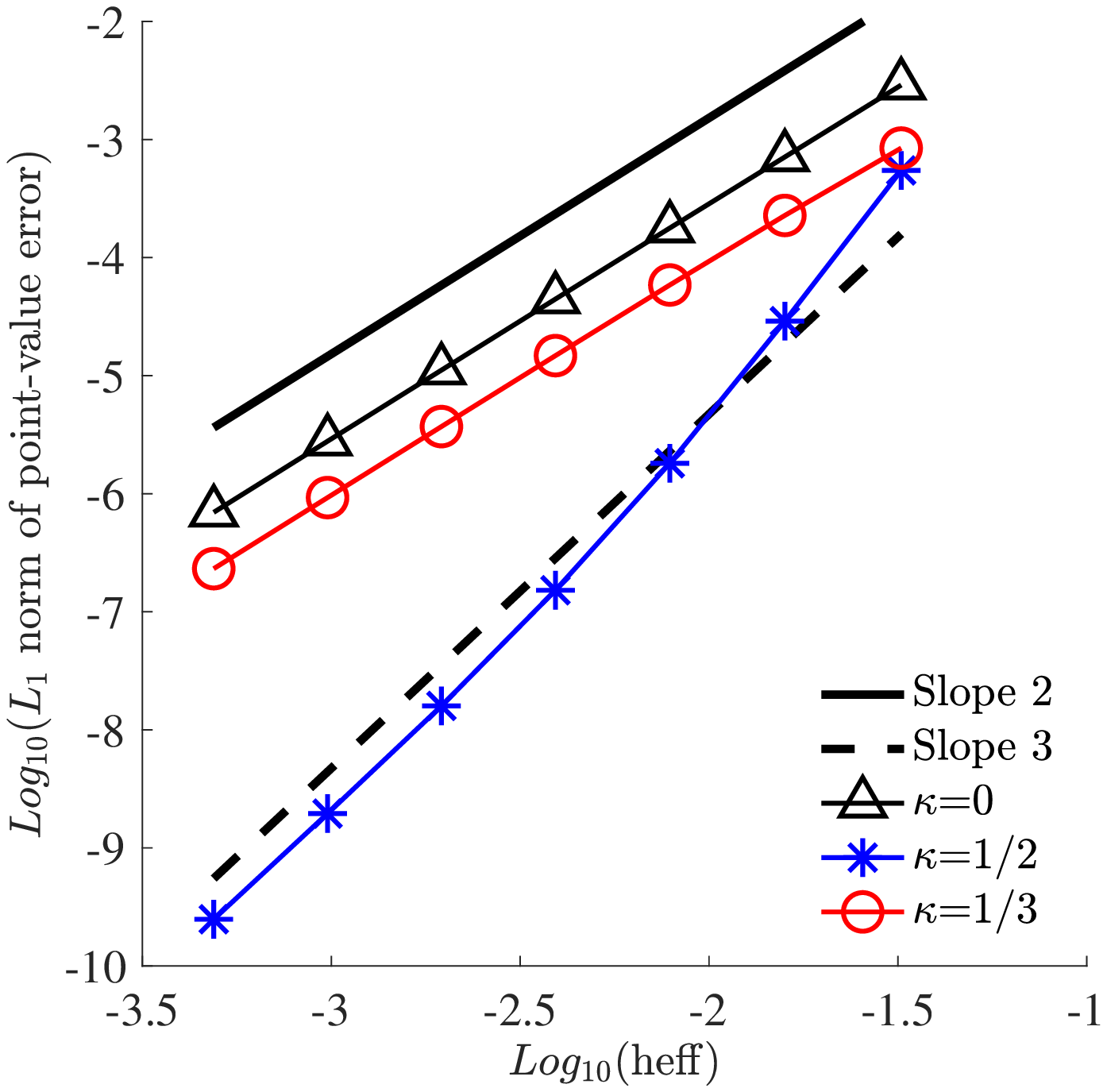}
          \caption{$L_1( {\cal E}_p ) $.}
          \label{fig:unsteady_err_p}
      \end{subfigure}
      \hfill
          \begin{subfigure}[t]{0.3\textwidth}
        \includegraphics[width=\textwidth]{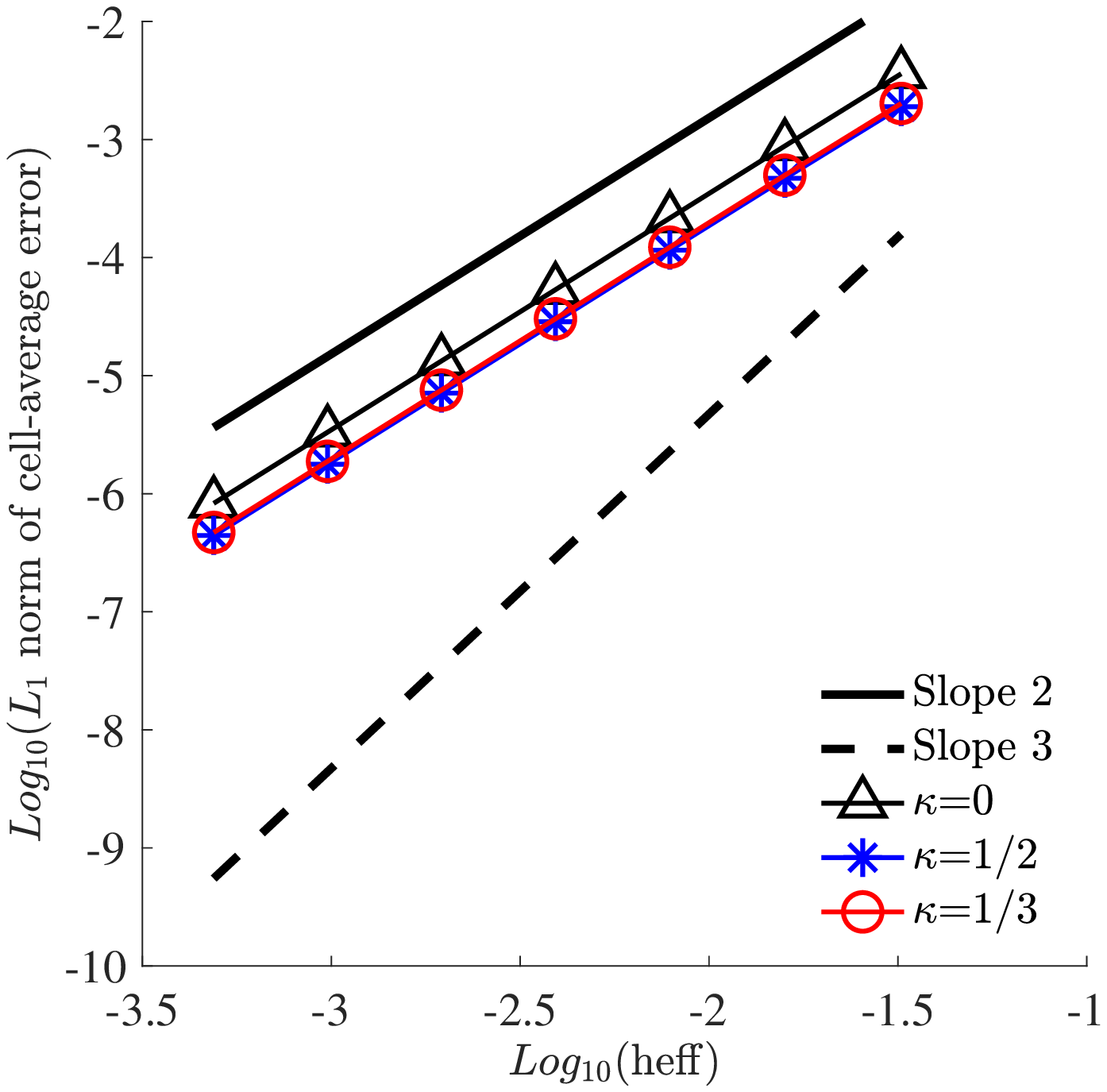}
          \caption{$L_1( {\cal E}_c) $.}
       \label{fig:unsteady_err_ca}
      \end{subfigure}
     \hfill
          \begin{subfigure}[t]{0.3\textwidth}
        \includegraphics[width=\textwidth]{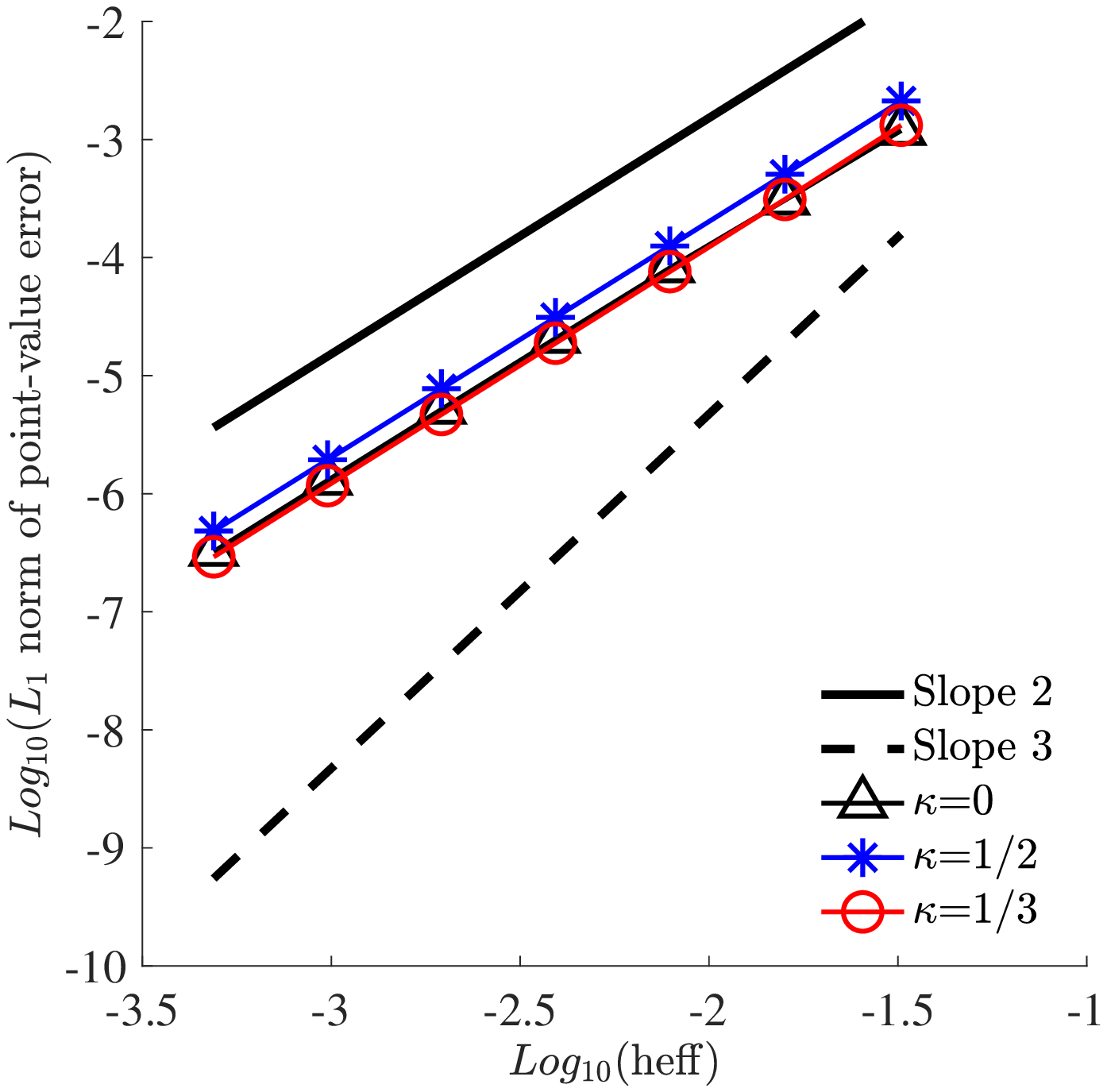}
          \caption{$L_1( {\cal E}_p) $, lumped mass matrix.}
       \label{fig:unsteady_err_p_limped}
      \end{subfigure}
            \caption{
\label{fig:unsteady_error_p}%
Coupled QUICK ($\kappa=1/2$): error convergence results for the case of the unsteady Burgers equation.
} 
\end{figure}

  \begin{figure}[t]
    \centering
          \begin{subfigure}[t]{0.35\textwidth}
        \includegraphics[width=\textwidth]{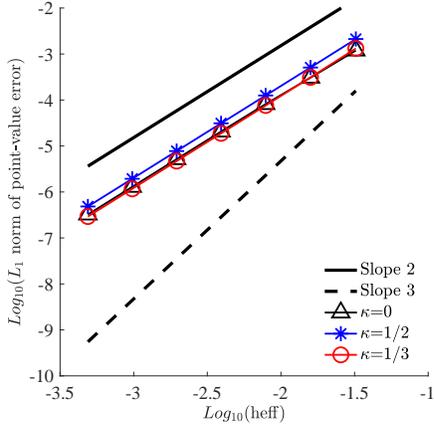}
          \caption{$L_1( {\cal E}_p) $: solution interpolation.}
       \label{fig:unsteady_err_ca_quickest_fd_sol}
      \end{subfigure}
      \hfill
          \begin{subfigure}[t]{0.35\textwidth}
        \includegraphics[width=\textwidth]{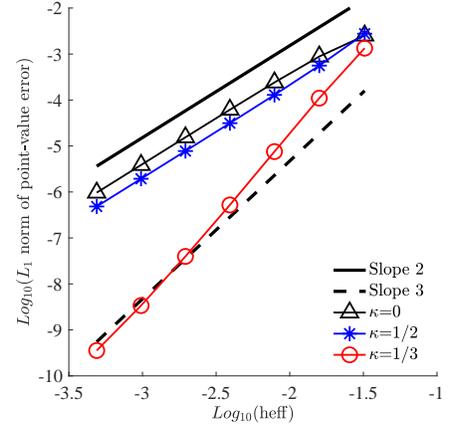}
          \caption{$L_1( {\cal E}_p ) $: flux interpolation.}
          \label{fig:unsteady_err_p_quickest_fd_flux}
      \end{subfigure}
            \caption{
\label{fig:unsteady_error_p_quickest_fd}%
QUICKEST ($\kappa=1/3$) applied in the finite-difference form (\ref{pointwise_scheme1aaaa}): error convergence results for the unsteady Burgers equation. Third-order only with the flux reconstruction ($\kappa=1/3$). 
} 
\end{figure}

  \begin{figure}[t]
    \centering
          \begin{subfigure}[t]{0.35\textwidth}
        \includegraphics[width=\textwidth]{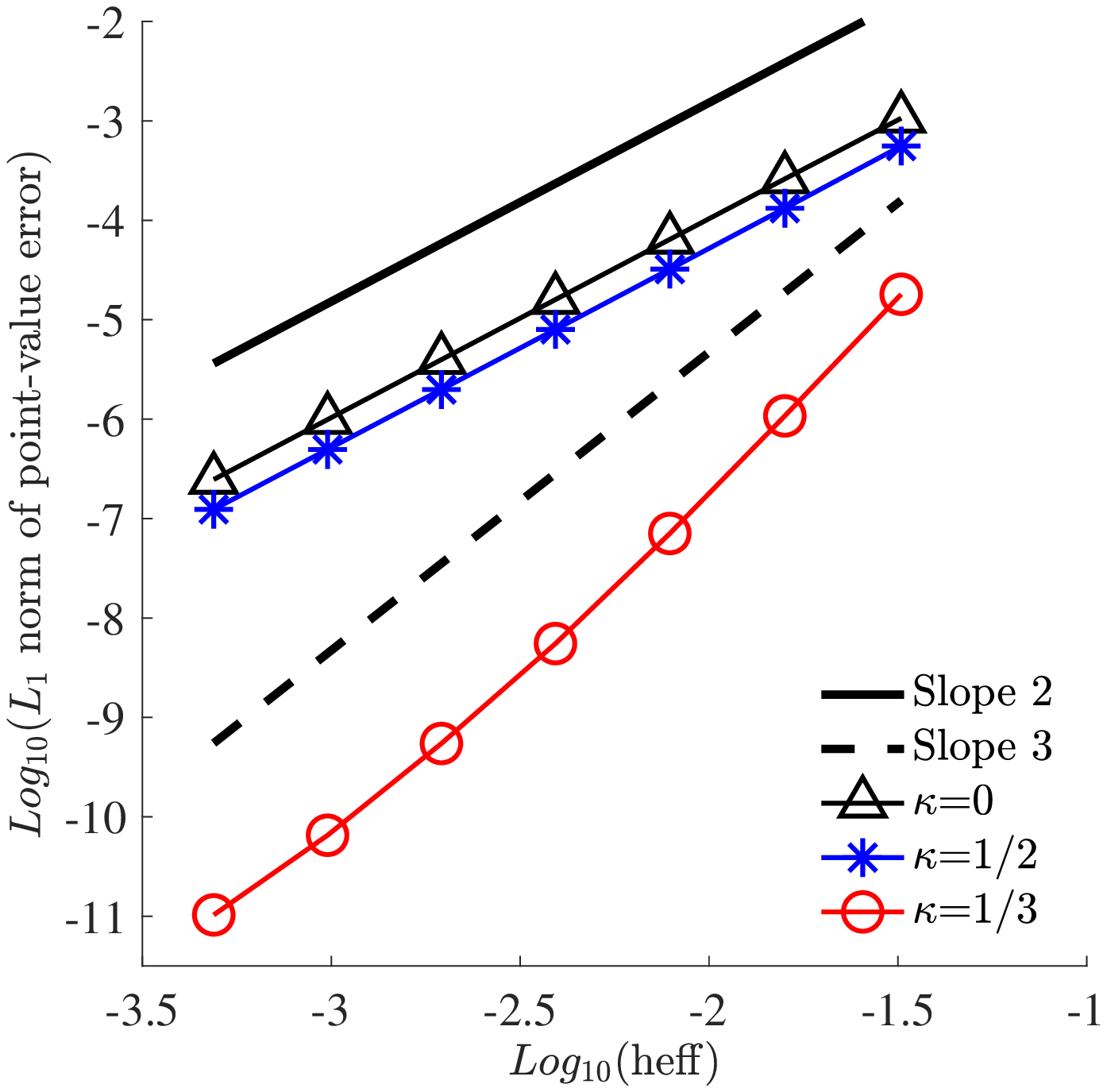}
          \caption{$L_1( {\cal E}_p) $: solution interpolation.}
       \label{fig:unsteady_err_ca_quickest_fd_sol_linear}
      \end{subfigure}
      \hfill
          \begin{subfigure}[t]{0.35\textwidth}
        \includegraphics[width=\textwidth]{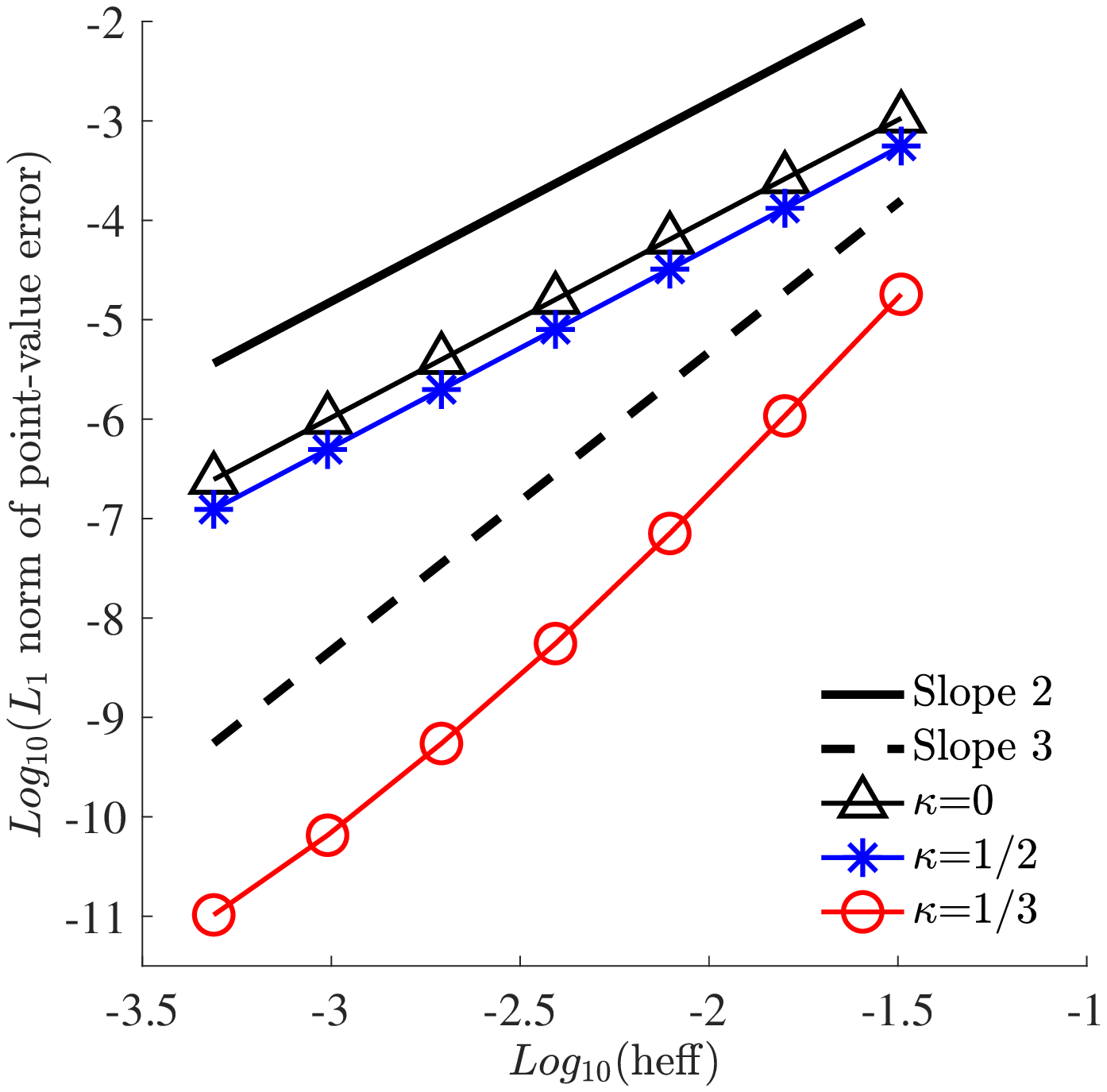}
          \caption{$L_1( {\cal E}_p ) $: flux interpolation.}
          \label{fig:unsteady_err_p_quickest_fd_flux_linear}
      \end{subfigure}
            \caption{
\label{fig:unsteady_error_p_quickest_fd_linear}%
QUICKEST ($\kappa=1/3$) applied in the finite-difference form (\ref{pointwise_scheme1aaaa}): error convergence results for the unsteady linear convection equation $u_t + (0.75 u)_x =  0$. Third-order with both the solution and flux reconstructions ($\kappa=1/3$). 
} 
\end{figure}

Figure \ref{fig:unsteady_error_p} shows results for the coupled QUICK scheme. 
As expected, third-order accuracy is achieved only with $\kappa=1/2$ and only in the point-valued solution: 
third-order convergence is observed in the point-valued error in Figure \ref{fig:unsteady_err_p}; 
second-order convergence in the cell-averaged error in Figure \ref{fig:unsteady_err_ca}. 
To demonstrate that the time-derivative coupling is critically important for third-order accuracy, we performed the same computation
with the lumped mass-matrix approach as described in Section \ref{third_order_summary}. 
Results are shown in Figure \ref{fig:unsteady_err_p_limped}. Clearly, third-order accuracy is lost, as expected. 

Next, we consider the QUICKEST scheme in the form (\ref{pointwise_scheme1aaaa}) 
integrated in time with the three-stage SSP Runge-Kutta scheme \cite{SSP:SIAMReview2001}.
Here, we compare two approaches for the flux computation. One is the flux evaluation with 
the solution interpolation with $\kappa=1/3$, and the other is the direct flux interpolation 
with Equations (\ref{kappa_uL_simple_flux}) and (\ref{kappa_uR_simple_flux}). 
Results are shown in Figure \ref{fig:unsteady_error_p_quickest_fd}. 
As expected, third-order accuracy is not achieved with the solution interpolation as shown in Figure \ref{fig:unsteady_err_ca_quickest_fd_sol}.
On the other hand, third-order accuracy is achieved with the flux interpolation for the QUICKEST scheme ($\kappa=1/3$) 
as clearly shown in Figure \ref{fig:unsteady_err_p_quickest_fd_flux}. Third-order achieved with the flux interpolation implies that 
 the QUICKEST scheme is third-order for linear equations, where the flux interpolation is equivalent to the solution interpolation.
Just for completeness, we show results obtained for a linear convection equation with $f = 0.75 u$ with the same initial solution in the same periodic domain. See Figure \ref{fig:unsteady_error_p_quickest_fd_linear}. Third-order accuracy is obtained with both the solution and flux reconstructions. 
These results indicate that the accuracy verification for a linear equation is not sufficient and it must be performed also for nonlinear equations.

\begin{figure}[t]
  \begin{center}
    \centering
          \begin{subfigure}[t]{0.35\textwidth}
        \includegraphics[width=\textwidth]{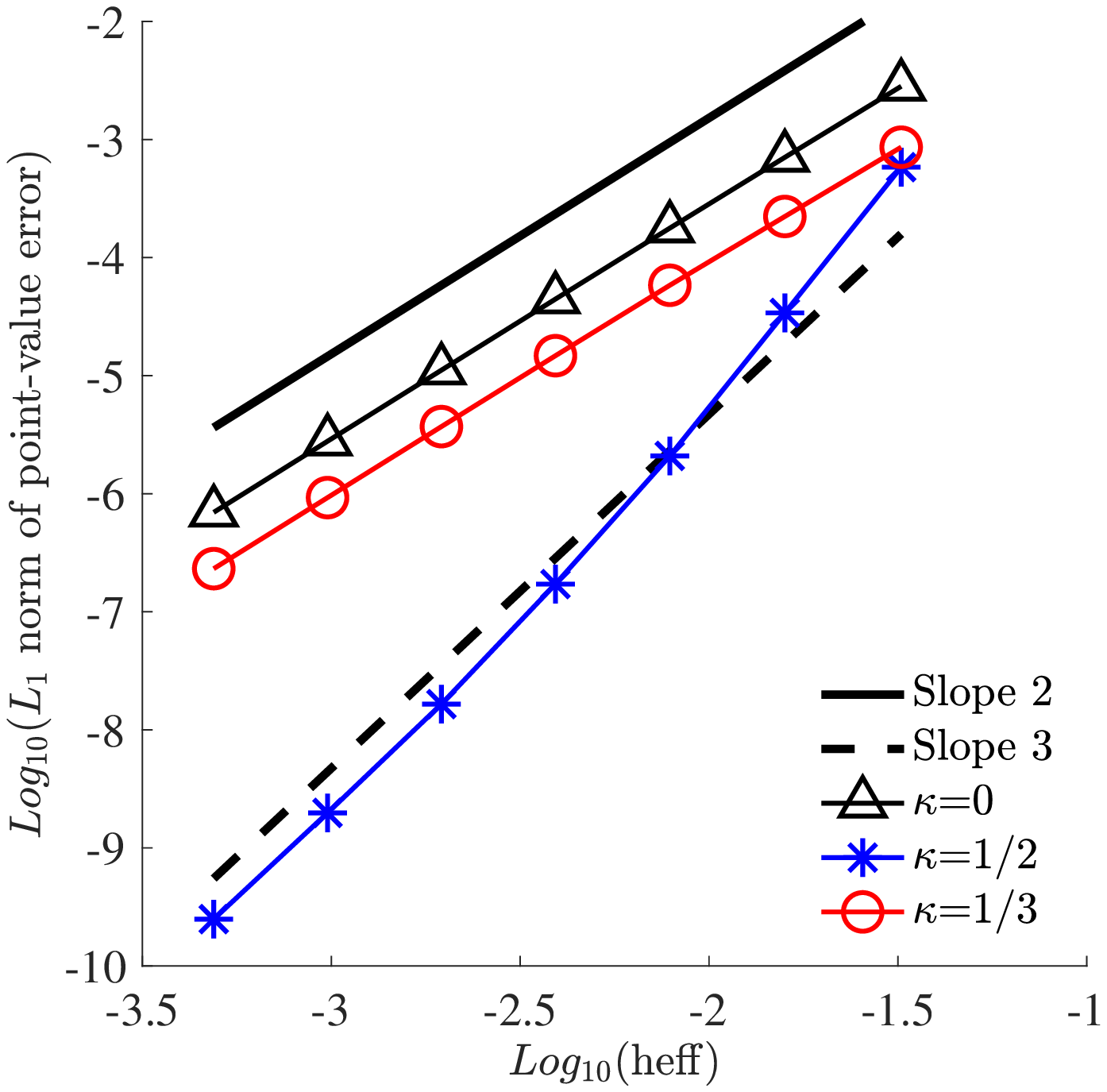}
          \caption{$L_1( {\cal E}_p ) $.}
          \label{fig:unsteady_err_p_quickest}
      \end{subfigure}
      \hfill
          \begin{subfigure}[t]{0.35\textwidth}
        \includegraphics[width=\textwidth]{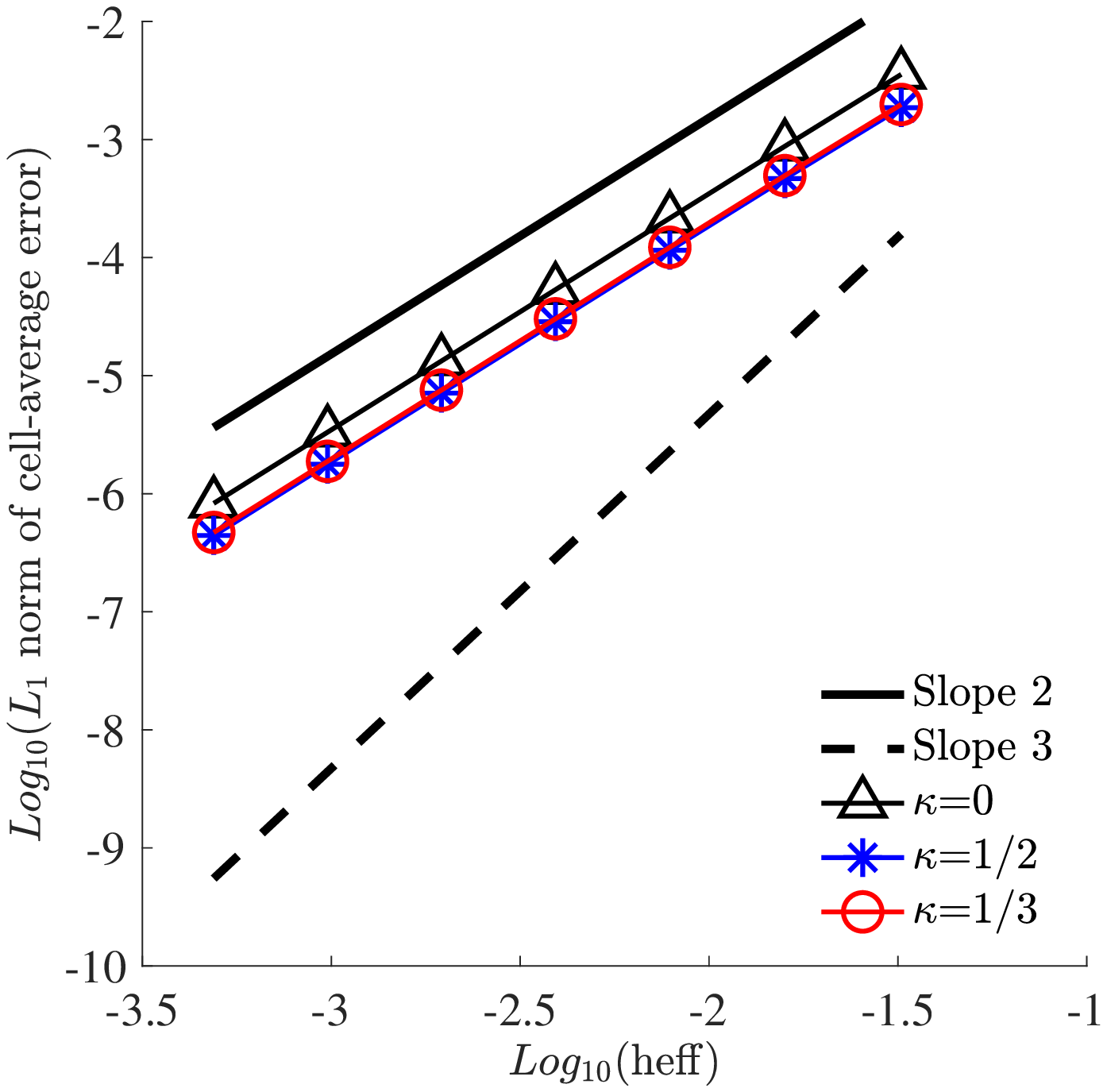}
          \caption{$L_1( {\cal E}_c) $.}
       \label{fig:unsteady_err_ca_quickest}
      \end{subfigure}
            \caption{
\label{fig:unsteady_error_p_quickest}%
Explicit QUICK of Van Leer ($\kappa=1/2$): error convergence results for the case of the unsteady Burgers equation.
} 
  \end{center}
\end{figure}

Finally, we tested Van Leer's explicit QUICK scheme mentioned in Section \ref{quickest_explained}, which is based on
the following approximation to Equation (\ref{pointwise_scheme0000}) without the forcing term $\overline{s}_i$:
\begin{eqnarray}
  \frac{d   {u}_i }{dt}   = -  \frac{1}{24} \left(
    \frac{d   \overline{u}_{i+1} }{dt}  - 2   \frac{d  \overline{u}_i }{dt} +   \frac{d   \overline{u}_{i-1} }{dt} 
   \right)   - \frac{1}{h}  [ F_{i+1/2}  - F_{i-1/2}  ]  ,
      \label{pointwise_scheme011}
\end{eqnarray}
where the time derivatives on the right hand side are evaluated explicitly as
\begin{eqnarray}
  \frac{d  \overline{u}_{i-1} }{dt} = - \frac{1}{h}  \left[  F_{i-1/2} -F_{i-3/2}  \right] ,   
    \frac{d  \overline{u}_{i} }{dt}   =   - \frac{1}{h}  \left[  F_{i+1/2} -F_{i-1/2}  \right] ,  
      \frac{d  \overline{u}_{i+1} }{dt}  = -   \frac{1}{h}  \left[  F_{i+3/2} -F_{i+1/2}  \right] ,
  \label{quickest_03}
\end{eqnarray}
which can be easily implemented in two steps: in the first step, we compute and store the residuals: e.g., in the cell $i$, 
\begin{eqnarray}
  Res_i^{(1)} =  \frac{1}{h}  \left[  F_{i+1/2} -F_{i-1/2}  \right] ,
\end{eqnarray}
where the superscript merely indicates the first step, then in the second step, we compute a corrected residual $Res_i$ 
to form the semi-discrete equation,
\begin{eqnarray}
  \frac{d   {u}_i }{dt}    +   Res_i = 0, 
  \quad  Res_i =   Res^{(1)}_i  - \frac{1}{24} \left(  Res^{(1)}_{i+1}  -2 Res^{(1)}_{i}+ Res^{(1)}_{i-1}    \right),
      \label{pointwise_quickest_00}
\end{eqnarray}
and integrate it in time with the three-stage SSP Runge-Kutta scheme.  Results are shown in Figure \ref{fig:unsteady_error_p_quickest}.
As expected, third-order accuracy is achieved only with $\kappa=1/2$ and only as a point-valued solution. 
See Figure \ref{fig:unsteady_err_p_quickest}. Third-order accuracy is not achieved as a cell-averaged solution as expected, 
which is shown in Figure \ref{fig:unsteady_err_ca_quickest}.

\section{Conclusions}
\label{conclusions}
 
We have clarified third-order accuracy of the QUICK scheme, which is a third-order finite-volume discretization of
the integral form of a conservation law with point-valued solutions stored as numerical solutions at cell centers.
Third-order accuracy in the point-valued solution has been verified by a detailed truncation error analysis for a general nonlinear conservation law and
also by a series of thorough numerical experiments. For unsteady problems, the QUICK scheme requires a consistent spatial discretization of the time derivative.
Two approaches have been considered: the coupled QUICK scheme and the explicit QUICKEST scheme of Leonard \cite{Leonard_QUICK_CMAME1979}.
The coupled QUICK scheme, a simple example of the deconvolution approach \cite{Denaro:IJNMF1996}, has been proved to be third-order accurate by a truncation error analysis and also demonstrated numerically for 
a nonlinear problem. On the other hand, the QUICKEST scheme has been shown to be a third-order conservative finite-difference scheme. 
As such, it is third-order accurate for linear equations (on uniform grids) but second-order accurate for nonlinear equations unless a direct flux reconstruction is performed. 
These facts have been shown by a truncation error analysis and also by numerical experiments. 
Also, we discussed confusions that can arise when the QUICK scheme is considered as a finite-difference scheme.
The QUICK scheme should not be interpreted as a finite-difference scheme (not proposed as such by Leonard), which can
cause lots of confusions, with the exception of the QUICKEST scheme. 
Finally, we briefly described an explicit QUICK scheme of Van Leer and presented numerical results to confirm third-order accuracy. 

Although not discussed, the QUICK scheme can be made to preserve third-order accuracy on irregular grids by fitting a quadratic polynomial over a set of 
irregularly-located solutions: $u_{i-1}$, $u_i$, and $u_{i+1}$ with $x_{i+1}-x_i \ne x_i - x_{i-1}$. Also, the coupled QUICK scheme can be extended to irregular grids by 
approximating the curvature term with the second derivative of the same quadratic polynomial. However, third-order accuracy is achieved only for the convective term with 
the quadratic interpolation. For diffusion, a cubic interpolation is required on irregular grids. 
Furthermore, the QUICKEST scheme is no longer third-order accurate on 
irregular grids even with the flux reconstruction \cite{Shu_WENO_SiamReview:2009,Merryman:JSC2003}. A further discussion 
on irregular grids will be provided elsewhere. 

This paper has shown together with the previous paper \cite{Nishikawa_3rdMUSCL:2020} that the a finite-volume discretization of the integral form is third-order accurate with $\kappa=1/3$ if cell-averaged solutions are stored as numerical solutions at cell centers (the third-order MUSCL scheme), and third-order accurate with $\kappa=1/2$ if point-valued solutions are stored as numerical solutions at cell centers (the QUICK scheme). In a subsequent paper, we will clarify fake third-order error convergence reported 
in the literature for he U-MUSCL scheme.

\addcontentsline{toc}{section}{Acknowledgments}
\section*{Acknowledgments}

The author is grateful to Emeritus Professor Bram van Leer for illuminating discussions and helpful suggestions,
and also to Professor Filippo Maria Denaro for clarifying the deconvolution approach. 
The author gratefully acknowledges support from Software CRADLE, part of Hexagon, the U.S. Army Research Office 
under the contract/grant number W911NF-19-1-0429 with Dr. Matthew Munson as the program manager, and 
 the Hypersonic Technology Project, through the Hypersonic Airbreathing Propulsion Branch of the NASA Langley
 Research Center, under Contract No. 80LARC17C0004.

\addcontentsline{toc}{section}{References}
\bibliography{../../bibtex_nishikawa_database}
\bibliographystyle{aiaa}

 
\end{document}